\documentclass[10pt]{amsart}
\usepackage{amscd,amssymb,graphics}
\usepackage[hidelinks]{hyperref}

\input xy
\xyoption{all}
\usepackage{epsfig}

\newtheorem{thm}{Theorem}[section]
\newtheorem{corol}[thm]{Corollary}
\newtheorem{lemma}[thm]{Lemma}
\newtheorem{prop}[thm]{Proposition}

\newtheorem{f}[thm]{Fact}

\theoremstyle{definition}
\newtheorem{defin}[thm]{Definition}

\theoremstyle{remark}
\newtheorem{remark}[thm]{Remark}
\newtheorem{remarks}[thm]{Remarks}
\newtheorem{example}[thm]{Example}
\newtheorem{examples}[thm]{Examples}
\numberwithin{equation}{section}

\newcommand{\delete}[1]{}    % Comment out text.

\def\e{{\varepsilon}}
\def\a{\alpha}

\newcommand{\nt}{\noindent}

\newcommand{\tv}{\tilde v}

\newcommand{\ri}{\rightarrow}

\newcommand{\bit}{\begin{itemize}}
\newcommand{\eit}{\end{itemize}}
\newcommand{\ben}{\begin{enumerate}}
\newcommand{\een}{\end{enumerate}}

\def\Homeo{{\mathrm{Homeo}}\,}

\def\R {{\mathbb R}}

\def\N{{\mathbb N}}
\def\Z {{\mathbb Z}}
\def\T{{\mathbb T}}

\def\Is{\operatorname{Is}}

\def\Asp{\operatorname{Asp}}
\def\RUC{\operatorname{RUC}}
\def\WRUC{\operatorname{WRUC}}

\def\LUC{\operatorname{LUC}}

\def\WAP{\operatorname{WAP}}
\def\AP{\operatorname{AP}}
\begin{document}

%%%\delete{
%%%%%%%%%%%%%%%%%%%%%%%%%%%%%%%%%%%%%%%%%%%%%%%%%%%%%%%%%%%%%%
%% This a placeholder for the TOPLOGY PROCEEDINGS logo %%%%%%
%\noindent                                             %%%%%%
%\begin{picture}(150,36)                               %%%%%%
%\put(5,20){\textbf{Topology Proceedings,}}                       %%%%%%
%\put(5,7){\textbf{145, 2004, 119-145.}}%%%%%%
%\put(0,0){\framebox(140,34)}                          %%%%%%
%\put(2,2){\framebox(136,30)}                          %%%%%%
%\end{picture}                                         %%%%%%
%%%%%%%%%%%%%%%%%%%%%%%%%%%%%%%%%%%%%%%%%%%%%%%%%%%%%%%%%%%%%
%%%%
%%%}

\title[]{\bf Fragmentability and representations of flows}
\author[]{Michael Megrelishvili}
\address{Bar-Ilan University, Israel}
\email{megereli@math.biu.ac.il}
\urladdr{http://www.math.biu.ac.il/$^\sim$megereli}
%\date{December, 2003}
%\date{May 29, 2005}
\date{May 29, 2005 (revised: September, 2019)}


%15June2006 Fact 5.2       S \to

\begin{abstract}
Our aim is to study weak$^*$ continuous representations of
semigroup actions into the duals of ``good'' (e.g., reflexive and
Asplund) Banach spaces. This approach leads to flow analogs of
Eberlein and Radon-Nikod\'ym compacta and a new class of functions
({\it Asplund functions}) which intimately is connected with
 representations on Asplund spaces and includes the class of weakly almost
periodic functions.
 We show that a flow is weakly
almost periodic iff it admits sufficiently many reflexive
representations.
 One of the main technical tools
in this paper is the concept of {\it fragmentability}
(which actually comes from Namioka and Phelps) and widespreadly used in
topological aspects of Banach space theory.
 We explore fragmentability
as ``a generalized equicontinuity'' of flows.
This unified approach allows us to obtain several
dynamical applications.
We generalize and strengthen some
results of Akin-Auslander-Berg, Shtern,
Veech-Troallic-Auslander and Hansel-Troallic. We establish that
frequently, for linear G-actions, weak and strong topologies
coincide on, not necessarily closed, G-minimal subsets. For
instance such actions are ``orbitwise Kadec``.
\end{abstract}

%\thanks{{\it 2000 Mathematical Subject Classification.}
%Primary 43A60, Secondary 22A05, 22C05, 22F05, 43A07.}

\maketitle

\tableofcontents

%\begin{equation} \label{f1}
%A+t=b
%\end{equation}

\section{Introduction}

 Every compact jointly continuous $G$-flow $X$ admits a faithful weak$^*$
 continuous Banach representation. More precisely, $X$ is
$G$-embedded into the dual ball $B(V^*)$ as a weak$^*$
compact $G$-subset of some Banach space $V$, where the group
$G$ acts continuously on $V$
by linear isometries. Indeed, this is a standard fact (see
Teleman's paper \cite{Te}, or for a more detailed discussion, the
survey \cite{Pe1}) for $V=C(X)$, where one can identify $x\in X$
with the point mass $\delta_x \in C(X)^*$. The geometry of $C(X)$,
in general, is bad. For example, a very typical disadvantage here is
the norm discontinuity of the dual action of $G$
on $C(X)^*$. One of the results of \cite{Mefr} guarantees
(see also Corollary \ref{c-dual} below) the norm
continuity of the dual action of the group $G$ on $V^*$ provided
that $V$ is Asplund. Recall that a Banach space $V$ is {\it
Asplund} iff the dual $A^*$ is separable for every separable
Banach subspace $A$ of $V$.

The following general question arises: how good can a Banach
space $V$ be among all possible w$^*$-continuous faithful
$G$-linearizations of $X$ into $V^*$ ? For instance when can $V$ be
chosen Asplund or reflexive ? We show that the reflexive case (for
second countable $X$) can be reduced completely to the question if
$X$ is a weakly almost periodic (in short: {\it wap}) flow.
%in the sense of Ellis-Nerurkar \cite{EN}.

{\it Eberlein compact} in the sense of Amir and Lindenstrauss \cite{AL} is
a compact space which can be embedded into $(V, weak)$ for some
Banach space $V$. It is well known \cite{DFJP}
that a compact space $X$ is
Eberlein iff it can be embedded into the unit ball $(B(V), weak)$
 of some reflexive space $V$. If $X$ is a weak$^*$ compact subset
in the dual $V^*$ of an Asplund space $V$ then, following Namioka
\cite{Na2},  $X$ is called {\it Radon-Nikod\'ym compact} (in
short: RN). Every reflexive Banach space is Asplund. Hence, every
Eberlein compact is RN.

Now introduce map versions of these concepts. Let $f: X \to X$ be
a selfmap on a compact space $X$. Let us say that $f$ is an {\it
Eberlein (Radon-Nikod\'ym) map} if it admits a weak$^*$
linearization into certain reflexive (resp.: Asplund) Banach
space. That is, there exists a reflexive (Asplund) Banach space
$V$ and a weak$^*$ embedding $X \hookrightarrow B(V^*)$ in such a
way that $f:X \to X$ is a restriction of the adjoint $F^*: V^* \to
V^*$ of some linear operator $F: V\to V$ which is {\it
non-expansive} ( $\|F\| \leq 1$). In this point of view, the space
$X$ is Eberlein or RN iff the identity mapping $1_X: X \to X$ is
Eberlein or RN, respectively.

Clearly, every metric compact space is Eberlein
%and even {\it uniformly Eberlein}
since it is a compact subset of the Hilbert space $l_2$. In
contrast, even simple maps on metric compacta can be
non-Eberlein. For example, the $f(x)=x^2$ map on the closed
interval $[0,1]$ is not Eberlein. The map
$$
f: \T^2 \to \T^2,
\hskip 0.2cm f([a],[b])=([a+b],[b+\sqrt{2}])
$$
defined on the torus $X = \T^2$
is not even RN (see Example \ref{skew}).

It is significant that a compact metric {\it cascade} $(\Z, X)$ is
wap (equivalently, is an Eberlein flow,
by virtue of Corollary \ref{seccwap}) iff the generating
selfhomeomorphism $f: X \to X$ leads to a {\it wap Markov
operator} $T_f: C(X) \to C(X)$ (see Downarowicz \cite{Do} and the
references there). The study of wap operators and corresponding
cascades goes back to the 60's
(K. Jacobs, B. Jamison, M. Rosenblatt, R. Sine, J. Montgomery, E.
Thomas and others).

The setting of maps and their linearizations admits a natural
generalization in terms of flow linearizations. We introduce below
 Eberlein and Radon-Nikod\'ym flows and show that a
 flow is weakly almost periodic
in the sense of Ellis-Nerurkar \cite{EN} iff it is a subdirect
product of Eberlein flows.
Investigation of RN flows naturally leads also to a new
class of functions which we call {\it Asplund functions}.
%It is a well known phenomenon that wap functions on
%flows and groups are hard to discover even for the group $\Z$ of
%all integers. The class of Asplund functions might provide a new
%useful tool .
%The paper is organized in the manner that
Our approach emphasizes more the
similarities (rather than the differences) between wap and Asplund
functions. We show that a function is wap (Asplund) iff it comes
from a matrix coefficient defined by a representation into
reflexive (resp.: Asplund) spaces. In both cases our
method is based on corresponding dualities and a factorization procedure
by Davis, Figiel, Johnson and Pelczynski
\cite{DFJP}. 
%160919 
In fact, in the present updated version we use its ``isometric version" 
by Lima, Nygaard and Oja  \cite{LNO}.

In the ``Asplund case'' the technical part uses
a modification
(using ``Asplund subsets'' instead of ``weakly compact'')
which is due to Stegall \cite{St}.

Let us briefly describe one of the ideas explored in the present paper.
Suppose that $X$ is a subflow (under some action by linear isometries)
of a weak$^*$ compact dual ball $B(V^*)$ of some Banach space $V$.
One of the important questions in Banach space theory is a relationship
between  norm and weak$^*$  topologies on $X \subset V^*$.
In the ``absolute case'' of the coincidence, we say that $X$
is a {\it Kadec subset} of $V^*$.
In such cases, $X$, as a flow, is equicontinuous. Conversely, every
compact metric
equicontinuous flow $X$ admits a flow representation in such a
way that $X$ becomes a norm compact subflow
(and hence a Kadec subset) of a suitable $B(V^*)$.
In general, as an attempt to
measure ``the level of equicontinuity'',
we can ask how close can two natural
topologies on $X$ inherited from $V^*$ be.
A more concrete and flexible enough question is:
for what flow representations is
the natural mapping $1_X: (X,weak^*) \to (X, norm)$ {\it fragmented}
in the sense of \cite{NP, JR}.
The latter means that every nonempty subset of $X$ admits
relatively weak$^*$ open nonempty subsets with arbitrarily
small diameters.

The great advantage of Asplund spaces is the (weak$^*$,
norm)-fragmentability of bounded subsets in their duals \cite{NP,
Na2}. Many modern investigations in Banach spaces concern Asplund
spaces, the notion of fragmentability and closely related
Radon-Nikod\'ym property (see \cite{Bour, Na2, DGZ, Fa, BL} and
the references therein). In \cite{Mefr,Meop} we study some
dynamical applications of fragmentability. In the present paper we
examine further developments exploring some ideas more familiar in
the topological aspects of Banach space theory.
% working in the unified framework of flow
%representations.

For the convenience of the reader we have tried to make the exposition
self-contained.

\delete{
Our exposition has been influenced mainly
by Namioka \cite{Na1}, Ellis and Nerurkar \cite{EN},
Junghenn \cite{Ju} and Fabian \cite{Fa}.
}

\vskip 0.1cm

\nt {\bf Acknowledgments}: I would like to express my gratitude
to M. Fabian, E. Glasner, A. Leiderman, V. Pestov and V. Uspenskij
for helpful comments and suggestions. The
main results of this paper were presented
at the 9th Prague Topological Symposium
(August, 2001) and also at the Auckland ``Summer''
Topological Conference (July, 2002).
I would like to thank the organizers for their kind
invitation and hospitality.
This work is supported by ISF grant no. 4699.

\section{Preliminaries}

\vskip 0.1cm

 The closure and the interior operators in topological spaces
will be denoted by $cl$ and $int$, respectively.
 If $A$ is a subset in a Banach space then $sp(A)$ is the linear span of $A$.
Let $\mu$ be a uniform structure on a set $X$. Its induced
topology on $X$ will be denoted by $top(\mu)$. A uniformity
%(we do not require the separation axiom)
$\mu$ on a topological space $(X,\tau)$ is said to be {\it
compatible} if $top(\mu)=\tau$.
A (left) {\it flow} $(S,X)$ consists of a topologized
semigroup $S$ and a (left) action $\pi : S\times X \to X$ on
a topological space $X$. We reserve the symbol $G$ for the case
when $S$ is a group. As usual we write simply $sx$ instead of
$\pi(s,x)={\tilde s} (x)= {\tilde x}(s)$. ``Action'' means that
always $s_1(s_2x)=(s_1s_2)x$. If $S$ is a monoid, we assume that
the identity $e$ of $S$ acts as the identity transformation of
$X$.
%We always require that the $s$-translation ${\tilde s}: X \to
%X, x \mapsto sx$ is continuous.
Every $x \in X$ defines an orbit
map ${\tilde x}: S \to X, s \mapsto sx$. Say that a topologized
semigroup $S$ is: (a) left (right) topological; (b)
semitopological; (c) topological if the multiplication function
$S\times S \to S$ is left (right) continuous, separately
continuous, or jointly continuous, respectively. Let $S$ be a
semitopological semigroup. A left flow $(S,X)$ is said to be a
{\it semitopological flow} if the action
is separately continuous.
%In jointly continuous case sometimes we say {\it $S$-space.

 A {\it right flow} $(X,S)$ can be defined analogously. If
$S^{opp}$ is the {\it opposite semigroup} of $S$ with the same
topology then $(X,S)$ can be treated as a left flow $(S^{opp},X)$
(and vice versa).

{\it If not stated otherwise the flows below are assumed to be
semitopological. ``Compact'' will mean compact and Hausdorff.}

%delete{
%(a) The semigroup $X^X$ is a right topological (compact)
%semigroup for every (compact) space $X$.
%
%(b) $(C_p(X,X), X)$ is a semitopological flow for every $X$.
%
%(c) $(C_p(X),S), \hskip 0.4cm (fs)(x)=f(sx)$
%is a semitopological (right) flow
%
%for every semitopological flow $(S,X)$.
%(d) $(\Theta(V),B(V))$ is a semitopological flow for every Banach space
%}

Let $h: S_1 \to S_2$ be a semigroup homomorphism, $S_1$ act
on $X_1$ and $S_2$ on $X_2$. A map $f: X_1 \to X_2$ is said to be
$h$-{\it equivariant} if $f(sx)=h(s)f(x)$ for every $(s,x) \in S_1
\times X_1$. For $S_1=S_2$ with $h=1_S$, we say {\it $S$-map}.
The map $h: S_1 \to S_2$ is an {\it antihomomorphism } iff
$h: S_1 \to S_2^{opp}$ (the
same assignement) is a homomorphism.

An  {\it $S$-compactification} of $(S,X)$ is a continuous $S$-map
$\alpha: X \to Y$ with a dense range ( {\it $S$-compactification
map}) into a compact $S$-flow Y. A (jointly continuous) flow
$(S,X)$ is said to be (resp.:{\it  joint continuously}) {\it
compactifiable} if there exists an $S$-compactification $\gamma: X
\to Y$ into a (jointly continuous) $S$-flow $Y$ such that $\gamma$
is a topological embedding. Following Junghenn
\cite{Ju} we define a {\it bicompactification} $m=(h,\alpha):
(S_1,X_1) \rightrightarrows (S_2,X_2)$ as a pair $h: S_1 \to S_2,
\alpha: X_1 \to X_2$, where $(S_2,X_2)$ is a semitopological flow
with compact $S_2$ and $X_2$, the map $h$ is a continuous
homomorphism and $\alpha$ is a continuous $h$-equivariant map with
a dense range.

Let $V$ be a Banach space with the dual $V^*$. Set
$$B(V)=\{v \in V:  ||v|| \leq 1 \} \hskip 0.1cm \text{and} \hskip 0.1cm
 %$$B(V^*)=\{v \in V^*:  ||v|| \leq 1 \} $$
\Theta(V)=\{\sigma \in L(V,V):  ||\sigma|| \leq 1 \}.$$
In most cases we endow the sets $B(V)$, $B^*=B(V^*)$
and $\Theta(V)$ with weak, weak$^*$ and weak operator topologies,
respectively.  Sometimes we use the subscripts ``w'' and
``w$^*$''. The subscript ``s'' will mean the strong operator
topology. The pairs $(\Theta(V)_w,B(V)_w)$ and
$(B(V^*)_{w^*},\Theta(V)_w)$ are semitopological flows.
The $\Theta(V)_s^{opp}$-flow
%$(\Theta(V)_s^{opp},
$B(V^*)_{w^*}$ is jointly continuous as it
follows directly from Fact \ref{Teleman}. It induces a right
action of the isometry group $\Is(V)=\{g \in Aut(V):  ||g|| = 1 \}$.
Alternatively, we have a left
action defined by
$$\Is(V) \times B^* \to B^*, \quad (gf)(v)=f(g^{-1}v)$$
Hence, $(\Is(V)_s, B^*_{w^*})$ is a well defined jointly continuous
action (see also Remark \ref{rem-anti}).

%\delete{
%\begin{lemma} \label{distel}
%(i) Let $X$ be a compact space. and
%let $K_1$ and $K_2$ be weakly compact subsets in C(X). Then
%$K_1\cdot K_2:=\{f_1f_2 | f_1 \in K_1, f_2 \in K_2 \}$
%is also weakly compact.
%
%(ii) If $K$ be a weakly compact subset and $||v|| \leq r_0 < 1$ for every
%$v \in K$ then $\cup \{ K^n | n\in \N \} \cup \{0\}$ is w-compact, where
%$K^n: = K\cdot K \cdots K$.
%\end{lemma}
%\begin{proof}
%(i) $K_1\cdot K_2$ is bounded.
%Since $K_1$ and $K_2$ are sequentially w-compact, it is clear that
% $K_1\cdot K_2$ is sequentially p-compact. Therefore, also
%w-compact (again by Eberlein-Smulian theorem).
%Thus we can apply Fact \ref{Groth}.
%
%(ii) easily follows from (i).
%\end{proof}
%}

The Banach algebra of all continuous real valued  bounded functions
on a topological space $X$ will be denoted by $C(X)$. The same set
with the pointwise topology ({\it p-topology})
is denoted by $C_p(X)$. Let $X$ be a
(left) $S$-flow then it induces the antihomomorphism $h: S \to
C(X)$ and the corresponding (right) action $C(X)\times S \to C(X)$
where $(fs)(x)=f(sx)$. In the case of a topological group $S=G$,
%%%(or, more generally, of a semigroup $S$ with a continuous involution)
we can define a homomorphism and a {\it left action} by
$(gf)(x)=f(g^{-1}x)$. While the translations are continuous, the
orbit maps ${\tilde f}: S \to C(X)$ are not necessarily (even
weakly) continuous. Denote by $\RUC_S(X)$ the set of all functions
$f \in C(X)$ such that the orbit map ${\tilde f}$ is norm
continuous. If we require only weak continuity, then we get the
definition of {\it weakly right uniformly continuous functions}
(see \cite{BJM}). Denote the corresponding set by $\WRUC_S(X)$.

The proof of the following fact is straightforward.

\begin{f} \label{f-c=r}
If $X$ is compact, then $\pi : S \times X \to X$ is jointly
continuous iff $C(X)=\RUC_S(X)$.
\end{f}

For general (separately continuous) action $\pi$, the set
$\RUC_S(X)$ is an $S$-invariant Banach subalgebra and the
corresponding Gelfand compactification $u_R: X \to X^R$ is a
universal (maximal) jointly continuous $S$-compactification of $X$.
If $X=S$ with the left regular action of $S$, then we simply
write $\RUC(S)$. If $S=G$ is a topological group, then $\RUC(G)$ is the
set of all usual right uniformly continuous functions. The algebra
of all left uniformly continuous functions (defined for the right
regular action of $S$ on $S$)) will be denoted by $\LUC(S)$.

The classical Gelfand-Naimark 1-1 correspondence between Banach
subalgebras of $C(X)$ and the compactifications of $X$ can be
extended to the category
 of jointly continuous $S$-flows
 using Banach $S$-subalgebras of $\RUC_S(X)$ (like the
 well-known results for topological group actions (see J. de Vries
 \cite{vr-embed})). One of the ways to verify this is to use
 the following fact which is a key idea of Teleman's above-mentioned
result, as well as in the paper of Uspenskij \cite{Usold}.

\begin{f} \label{Teleman}
Let $V$ be a Banach space. Suppose that a topologized semigroup
$S$ acts on $V$ from the right by linear non-expansive operators.
The following are equivalent:

(i) $V \times S \to V$  is norm jointly continuous.

%(ii) All orbit maps ${\tilde v}: S \to V$ are norm continuous;

(ii) The induced affine action $S \times B^*_{w^*} \to B^*_{w^*}$ is
jointly continuous.
\end{f}
\begin{proof} The dual action defines an injective antihomomorphism
of $\Theta(V)$ into $C(B^*,B^*)$. Now observe that the strong
operator topology on $\Theta(V)$ coincides with the compact open
topology inherited from $C(B^*,B^*)$.
\end{proof}

Recall the definition of weakly almost periodic functions and some
relevant facts.

\begin{defin} \label{wapdef}
Let $S$ be a semitopological semigroup and $X$ be an $S$-flow.
\bit \item [(i)] A function $f \in C(X)$ is said to be {\it weakly
almost periodic}, ({\it wap}, in short) if the orbit $fS=\{fs:
\hskip 0.1cm s \in S \}$ is relatively weakly compact in $C(X)$.
Write $f \in \WAP_S(X)$.

%(ii) $(S,X)$ is wap $\stackrel{def}{=}$ $C(X)=\WAP_S(X)$

% (iff $E(S,X) \subset C(X,X)$ (in part., $E(S,X)$ is semit.)).
\item
[(ii)] We say that $X$ is $S$-wap, or, $(S,X)$ is wap (notation:
$X \in [wap]^S$) if $\WAP_S(X)$ separates points and closed subsets
of $X$.
\item
[(iii)] We say that $S$ is {\it wap} (and write: $S \in [wap]$) if
the regular left action $(S,S)$ is wap.
\eit
\end{defin}

This general form of definition (i) can be found in the work of
Junghenn \cite{Ju}. For the left action $(S,S)$ we get the
classical notion of wap functions on $S$ (see Eberlein \cite{Eb}
and de-Leeuw Glicksberg \cite{LG1}). We use the notation $\WAP(S)$
instead of $\WAP_S(S)$.

Replacing 'weakly compact" in Definition \ref{wapdef} by "norm
compact" we get the definitions of {\it almost periodic} functions
and corresponding $S$-algebras $\AP_S(X), \AP(S)$.

%By Fact \ref{BJM},
%$WAP_S(X)$ is a Banach $S$-subspace in $C(X)$.
%Lemma \ref{distel} easily implies also that
%$WAP_S(X)$ is even an algebra $(f_1f_2)S \subset cl_w(f_1S)\cdot
%cl_w(f_2S)=K_1 \cdot K_2$).

%\begin{f} \label{simplest}
%For every
%\end{f}

Grothendieck's criteria \cite{Gr} for relative weak compactness
 leads to the following assertion.

\begin{f} \label{f-DLP} (Grothendieck's DLP)
A function $f \in C(X)$ defined on some $S$-flow is wap iff the following
{\it Double Limit Property} is satisfied:

\nt ($DLP$) \hskip 0.1cm For every pair of sequences $s_m \in S$
and $x_n \in X$
$$\lim_m \lim_n f(s_n x_m)=\lim_n \lim_m f(s_n x_m)$$
holds whenever both of the limits exist.
\end{f}

Recall also the following very useful fact.

\begin{f}
\label{Groth}
(Grothendieck's Lemma)
%\cite{gro})
Let $X$ be a compact space.
Then a bounded subset $A$ of $C(X)$ is
w-compact iff $A$ is $p$-compact.
\end{f}

The set $\WAP_S(X)$ is a Banach $S$-subalgebra in $C(X)$.
This is mentioned in \cite{Ju}.
The proof can be done using Fact
%\ref{f-DLP},
\ref{Groth} and the Eberlein-Smulian theorem.

%
%The less trivial is to show that $fh \in WAP_S(X)$ for every $f,h
%\in WAP_S(X)$. For compact $X$ this can be done by Fact \ref{Groth}, using the
%inclusion $(fh)S \subset (fS) \cdot (hS)$ and
%Eberlein-Smulian theorem. Since $WAP_S(X)$ is a
%subalgebra of $C(X)$,  compact $S$-flow $X$ is wap iff
%$WAP_S(X)$ separates the points of $X$ iff (use the
%Stone-Weierstrass Theorem) $C(X)=WAP_S(X)$.
%For general $X$ we can use Grothendieck's DLP.
This implies that our general Definition \ref{wapdef}(ii),
for compact $X$, is
equivalent to the definition of wap flows in the sense of
Ellis-Nerurkar \cite{EN}. Gelfand's compactification $u_W:~X \to
X^W$ induced by the algebra $\WAP_S(X)$ is the {\it universal wap
compactification} of $X$ (see for details \cite[Theorem 3.1]{Ju}).
 In particular, for the left regular action $(S,S)$ we get
the {\it universal wap semigroup compactification}
$u_W: S \to S^W$. It is
important to note that in this case $S^W$ is a compact
semitopological semigroup and enjoys the corresponding
universality property. By
our definitions, the flow $X$ (or, the semigroup $S$) is wap iff
$u_W$ is a topological embedding.

\vskip 0.1cm {\it Ellis semigroup} $E(S,X)$ (or, simply: $E(X)$)
for compact $X$ is the pointwise closure of the set of all
$s$-translations $\{ {\tilde s}: X \to X : \hskip 0.1cm s \in S
\}$ in the compact semigroup $X^X$. Denote by $\lambda: S \to
E(X), \lambda(s)={\tilde s}$ the corresponding natural
homomorphism. In general, $E(X)$ is only  {\it right topological},
that is, only the right translations $E(X) \to E(X), \hskip 0.1cm
s \mapsto sp$ are necessarily continuous.

\begin{f} \label{EN}
\cite{EN, Ju} For a compact $S$-flow $X$
the following are equivalent:

\ben \item $X$ is wap (that is, $C(X)=\WAP_S(X)$).
\item Each element of $E(S,X)$ is continuous ({\it
quasiequicontinuous} in terms of \cite{BJM}).

%(3) The semigroup $E(X)$ is semitopological ?????;
\item The pair $(E(X), X)$ is a semitopological flow with the compact
semitopological semigroup $E(X)$.

\item  There exists a bicompactication $(h,\alpha): (S,X)
\rightrightarrows (P,Y)$  such that  $\alpha: X \to Y$ is an
embedding and $h: S \to P$ is a semigroup compactification into a
compact semitopological semigroup $P$. \een
\end{f}
\begin{proof}
The principal implication (1) $\Rightarrow$ (2) directly follows
from Proposition \ref{ii} and Remarks \ref{r-two}(b).

For (4) $\Rightarrow$ (1) it suffices to show that
$C(Y)=\WAP_P(Y)$. Let $f \in C(Y)$ then the $P$-orbit $fP$ is
bounded. Since $P$ is compact then $fP$ is pointwise compact in
$C(Y)$. By Fact \ref{Groth}, $fP$ is even $w$-compact. Thus, $f
\in \WAP_P(Y)$.

Other implications are trivial.
\end{proof}
For compact $X$, Definition \ref{wapdef}(ii) agrees with the item (1)
in Fact \ref{EN}, as it easily follows by Stone-Weierstrass theorem.

If $X \in [wap]^S$ then $Y \in [wap]^P$
for every subsemigroup $P$ of $S$ and every
$P$-subflow $Y$ of $X$. Moreover, $[wap]^S$
is closed under {\it subdirect products} (subspaces of products).
The class of compact $S$-wap
flows is closed also under quotients.

\begin{f} \label{f-functions} \
\ben \item [(i)]
$\WAP_S(X) \subset \WRUC_S(X)$
 for every semitopological $S$-flow $X$.
Hence,  $\WAP(S) \subset \WRUC(S)$
for every semitopological semigroup $S$.
\item
[(ii)] If $G$ is a semitopological group then $\WAP_G(X) \subset
\RUC_G(X)$. In particular, $\WAP(G) \subset \RUC(G)$ holds.
\een
\end{f}
\begin{proof}
(i) The orbit map ${\tilde f}: S \to C(X)$ is clearly
p-continuous. If $f \in \WAP_S(X)$ then $cl_w(fS)$ is weakly
compact. Hence, $(fS,w)= (fS, p)$. Therefore ${\tilde f}$ is also
weakly continuous.

(ii) Follows from Fact 2.5 and Remark 2.6 (a) of \cite{Meop}.
It can be seen easily also as a corollary of Theorem
\ref{c-waplight} below.
\end{proof}
The inclusion $\WAP(G) \subset \RUC(G)$ is well known (see for
example \cite{He1} or \cite[Theorem 4.10]{BJM}).
Another proof of the inclusion
$\WAP_G(X) \subset \RUC_G(X)$ can be derived also by results of
\cite{JP}.

 \vskip 0.1cm
 For every reflexive Banach space
$V$ the semigroup $\Theta(V)$ is a weakly compact semitopological
semigroup \cite{LG1}. Observe that for every vector $v\in V$ with
norm 1, the orbit $\Theta(V)v$ of $v$ coincides with $B(V)$. This
guarantees the converse: if $\Theta(V)_w$ is compact then $B_w$ is
compact , and, hence, $V$ is necessarily reflexive.
%
%\delete{
%\begin{f} Let  $V$ be a Banach space . The following are equivalent:
%
%(1) $V$ is reflexive ;
%
%(2) $\Theta(V)_w$ is compact.
%
%(3) $B_V$ is $w$-compact.
%\end{f}
%}
For every reflexive $V$, the flows
$(\Theta(V), B(V))$ and  $(\Theta(V)^{opp}, B(V^*))$
are semitopological and (bi)compact. Hence, wap by Fact \ref{EN}.
%We call it a {\it canonical wap flow}.

\vskip 0.1cm

 One of our applications below (see section 8)
provides a simple  proof of
the following important theorem of Lawson \cite{La1} which in
itself is a generalization of Ellis theorem.

\begin{f} \label{lawson} (Ellis-Lawson's Joint Continuity Theorem).

Let $G$ be a subgroup of a  compact semitopological monoid $S$.
Suppose that $(S,X)$ is a semitopological flow with compact $X$.
Then the action $G\times X \to X$ is jointly continuous
and $G$ is a topological group.
\end{f}

A (not necessarily compact) $G$-flow $X$ is said to be {\it
minimal} if every orbit $Gx$ is dense in $X$. Equicontinuous
compact flows are the simplest one in Topological Dynamics. Every
equicontinuous compact flow is wap. The converse is true for every
minimal compact wap $G$-flow $X$ \cite{Tr, Au}. Below we show
(Theorem \ref{asplundr}) that
the compactness assumption is superfluous here. That is, every
minimal wap (and even, RN-approximable), not necessarily compact,
$G$-flow is equicontinuous.

%\delete{
%
%\nt {\bf Examples and properties:}
%
%\nt $\bullet$ wap flows:
%
%(a) $(P,Y)$ \quad where $P,Y$ are compact.  In particular,
%
%(b) Canonical wap: \quad $(\Theta(V),B_V)$
%For every {{\underline {reflexive}} Banach space $V$.
%
%%\ssk
%% $Is(V))=\{g \in Aut(V): ||g||=1 \}$
%
%%\nt (Grothendieck cr. $\Longrightarrow f|_B \in WAP_{\Theta}(B) \hssk
%%\forall f \in V^*$)
%\nt $\bullet$ Not wap flows:
%
%(a) $(\Z, \beta \Z);$
%
%Moreover, $(G, \beta_G G)$ for every non precompact $G$
%(follows from \cite{MPU}).
%
%(b) $(\Z, \{0,1\}^{\Z});$ \quad
%
%(c) $(\Z, [0,1]),  \hssk (\sigma(t):=t^2)$
%
%(d) $((0,\infty), \R)$
%
%(e) Every minimal but not equicontinuous $G$-flows.
%(For ex. , Furstenberg flow, nil-spaces).
%
%(f) Coset space $(G , G/H)$ where the action
%is not equicontinuous.
%
%(g) The equicontinuity is not suff. even for trivial $H$.
%Indeed, $(G,G)$ is not wap for $G=H_+[0,1]$ \cite{Merup}.
%}}
%

\section{Banach representations and matrix coefficients}
\label{s-repr}

\vskip 0.1cm
$\langle$ 
Let $V$ be a Banach space with the canonical duality $\langle ,\rangle: V
\times V^* \to \R$. If a semigroup $S$ acts from the right on $V$
(equivalently: if we have an antihomomorphism $S \to L(V,V)$) then
it induces a left action of $S$ on the dual $V^*$ such that
$\langle  vs, \psi \rangle=\langle v,s\psi \rangle$ for every $v \in V$ and $\psi \in V^*$.

\begin{defin} \label{def-repr}
A (non-expansive) {\it $V$-representation} of a flow $(S,X)$ is an
equivariant pair
$$(h, \alpha): (S,X) \rightrightarrows (\Theta(V)^{opp}, B^*)$$
%qmich  B^* or  V^* ?
where $h: S \to \Theta(V)^{opp}$ is a weak continuos homomorphism
(equivalently: antihomomorphism $S \to \Theta(V)$) and $\alpha:
X \to B^*$ is weak$^*$ continuous and equivariant, that is
$\alpha(sx)=h(s) \alpha(x)$. 
A little bit more flexible definition is replacing $B^*$ by $V^*$. 

We say that a representation is {\it strongly continuous} if $h: S
\to \Theta(V)_s$ is continuous. {\it Topologically faithful} (or,
simply: {\it faithful}) will mean that $\alpha: X \to (B^*,w^*)$
is a topological embedding.
\end{defin}

Let ${\mathcal K} \subset BAN$ be a subclass of Banach spaces.  
We say that a flow $(S,X)$ is:

\begin{itemize}
\item [(a)] ${\mathcal K}$-{\it representable}
if there exists a {\it faithful}
$V$-representation of $(S,X)$
%into $(\Theta(V)^*, B(V^*)_{w^*})$
for some $V \in {\mathcal K}$.
\item [(b)] ${\mathcal K}$-{\it approximable}
if there exists a  system $(h_i,{\alpha}_i)$
of representations of $(S,X)$ in $V_i$
%into $(\Theta(V_i), B(V^*_i)_{w^*})$)
separating points and closed subsets in $X$
with $V_i \in {\mathcal K}$
(equivalently, if $X$ is a subdirect product of
${\mathcal K}$-representable $S$-flows).
\item [(c)] {\it Eberlein} if it is REFL-representable.
\item [(d)] {\it Radon-Nikod\'ym} (in short: RN) if it is
ASP-representable.
\item [(e)]
RN-{\it approximable} if it is ASP-approximable.
\end{itemize}
%(e) $(S,X)$ is
%{\it uniformly Eberlein} if it is ${\mathcal Hilbert}$-representable.

In this definition REFL and ASP mean the classes of all reflexive
and Asplund spaces respectively. Since REFL $\subset$ ASP, every
Eberlein flow is RN. If $S$ is a trivial monoid and $X$ is compact,
then the definitions (c) and (d) give
exactly the classical notions of
Eberlein and RN compacta mentioned in the introduction.

\begin{remark} \label{r-automat}
 Sometimes weak continuous (anti)homomorphisms automatically
are strongly continuous. This happens for instance if either: (a)
$S$ is an arbitrary semitopological group and $V$ is reflexive;
(b) $S$ is a locally compact Hausdorff topological group; or (c)
$S$ is a topological group metrizable by a complete metric. The
first assertion follows from \cite[Theorem 2.8]{Meop} (or, from
Corollary \ref{PCP} below). For the last two assertions see
\cite{BJM}. For some other results of the nature ``weak implies
strong'' see also \cite{He1, He2, La2, Mefr, Meop}.
\end{remark}

%\nt $\bullet$ (v) $(S,X)$ is
%{\it residually finite dimensional} if
%in (ii) (resp., in (iii)) $V$ (resp., each $V_i$)
%is finite dimensional.

\vskip 0.1cm

%\delete{
%If $S=\{e\}$ (trivial monoid) then:
%$(\{e\},X)$ is $\{REFL\}$-approximable (or, wap) iff $X$ is Tychonoff;
%$(\{e\},X)$ is $\{REFL\}$-representable iff
%there exists a proper compactification $Y$ of $X$ such that $V$
%is an Eberlein compact.
%%Moreover, $(\{e\},X)$ is wap iff $X$ is Tychonov.
%
%%Therefore it is easy to find
%%Eberlein (and wap) which is not Eberlein.
%%Take for example,
%%$S=\{e\}, X=$nonmetrizable compact group, or $X=\beta \N$.
%
%If $X$ is compact then $X$ as a $\{e\}$-flow is $(RN)$ iff $X$ is $(RN)$
%in the sense of Namioka. This fact justifies our terminology.
%}

%\nt {\bf Remark:} about (wap functions):

%\nt $\bullet$  Ruppert:

%\nt $\bullet$ via Grothendieck's DLP (double limit property)

%\nt $\bullet$  ap functions on $S$: characters, ...

The following standard fact (see for
example \cite{Te}) states actually that every jointly
continuous action on compact spaces admits a faithful Banach
representation.

\begin{lemma} \label{genTel}
Let $(S,X)$ be a jointly continuous semigroup action on a compact
$X$. Then there exists a Banach space $V$ and a faithful strongly
continuous representation $(h, \alpha)$
of $(S,X)$ into the jointly continuous
affine action $(\Theta(V)_s^{opp}, B(V^*)_{w^*})$.
\end{lemma}
\begin{proof}
 Take $V=C(X)$ and define the antihomomorphism $h: S \to
\Theta(V)$ induced by the natural right action $C(X) \times S \to
C(X)$. This action is norm continuous by Fact \ref{f-c=r}
because $\RUC_S(X)=C(X)$.
Thus, $h$ is strongly continuous by Fact \ref{Teleman}.
Finally define the natural weak$^*$ embedding $\a: X \to
(B(C(X)^*)$ identifying each $x \in X$ with the point mass
$\delta_x \in B(C(X)^*)$.
\end{proof}

For every weakly continuous antihomomorphism $h : S \to L(V,V)$
and every chosen pair of vectors $v \in V$ and $\psi \in V^*$, there
exists a canonically associated {\it (generalized) matrix
coefficient}
$m_{v, \psi}: S \to \R, \hskip 0.3cm s \mapsto 
\langle vs, \psi \rangle=\langle v, s \psi \rangle$

\begin{equation*}
\xymatrix {
S \ar[d]_{h} \ar[r]^{m_{v,\psi}} & \R  \\
L(V,V) \ar[r]^{\tilde{v}} & V \ar[u]^{\psi}}
\end{equation*}

%nd
%\begin{diagram}[height=2.0em,width=3.0em,tight]
% &S &\rTo^{m_{v,\psi}} & \R  \\
% &\dTo^{h} &    &\uTo_{\psi}  \\
% & L(V,V)    & \rTo^{\tilde v} & V \\
%\end{diagram}

\begin{remark} \label{rem-anti}
In many important cases we can use homomorphisms instead of
antihomomorphisms. Indeed, if $S$ is a topological group (or, a
semigroup with a continuous involution), then we can define a {\it
homomorphism} $h^*: S \to L(V,V), s \mapsto h(s^{-1})$ and
redefine the function $m_{v,\psi}$ by $s \mapsto \langle s^{-1}v,\psi \rangle$.
\end{remark}

 It is natural to expect that matrix coefficients reflect good
properties of flow representations (see, for example, \cite{Pe1}).
We recall two well-known facts. The first example is the case of
Hilbert representations. If $h: G \to \Is(H)$ is a group
representation into Hilbert space $H$ and $\psi=v$, then the
corresponding map $g \mapsto \langle g^{-1}v, v \rangle$ is a {\it positive definite
function} on $G$. The converse is also true: every continuous
positive definite function comes from some continuous Hilbert
representation. Every positive definite function is wap (see
\cite{Bu}).

The second example comes from Eberlein \cite{Eb}
(see also \cite[Examples
1.2.f]{BJM}). If $V$ is reflexive, then every bounded
$V$-representation $(h, \alpha)$ and arbitrary pair $(v,\psi)$
lead to a weakly almost periodic function $m_{v,\psi}$ on $S$.
This follows easily by the (weak) continuity of the natural
operators defined by the following rule. For every given $\psi \in
V^*$ ($v \in  V$) define {\it introversion type operators} by
$$L_{\psi}: V \to C(S) \hskip 0.1cm \text{and} \hskip 0.1cm
R_v: V^* \to C(S), \text{where} \hskip 0.1cm  L_{\psi}(v)=
R_v(\psi)=m_{v,\psi}.$$

We say that a vector $v \in V$ is strongly (weakly) continuous if
the corresponding orbit map $\tv : S \to V, \tv(s)=vs$,
defined through $h: S \to \Theta(V)$, is strongly (weakly) continuous.

\begin{f} \label{introversion}
Let $h: S \to L(V,V)$ be a weakly continuous antihomomorphism with
the norm bounded range. Then

\ben
\item
$L_{\psi}: V \to C(S)$ (and $R_v: V^* \to C(S)$) are linear
bounded $S$-operators between right (left) $S$-actions.

\item
If $\psi$ (resp.: $v \in V$) is norm continuous, then $ m_{v,\psi}$
 is left (resp.: right) uniformly continuous on $S$.

\item
If $V$ is reflexive, then $m_{v,\psi} \in \WAP(S)$.
\een
\end{f}
\begin{proof} (1) Is straightforward.

(2) Since  $h(S)$ is norm bounded, 
$sup\{ ||vt||: \hskip 0.1cm t \in S \}=c < \infty$.
 Let $\psi$ be a norm continuous vector. In
order to establish that $m_{v,\psi} \in \LUC(S)$, observe that
$$|m_{v,\psi}(ts) - m_{v,\psi}(ts_0)|=| \langle vts, \psi \rangle - \langle vts_0,
\psi \rangle|=$$ 
$$ |\langle vt,s\psi \rangle - \langle vt,s_0\psi \rangle|
\leq ||vt|| \cdot 
||s\psi-s_0\psi|| \leq c \cdot ||s\psi-s_0\psi||.$$

Similar verification is valid for the second case.

(3) If the orbit $vS$ is relatively weakly compact in $V$ (e.g.,
if $V$ is reflexive), then the same is true for $L_{\psi}(vS)=
m_{v,\psi}S$ in $C(S)$. Thus $ m_{v,\psi}$ is wap. 
\end{proof}

\begin{f} \label{comesEb}
Let $(h, \alpha): (S,X) \rightrightarrows (\Theta(V)^{opp}, B^*)$ be an
equivariant pair with weak$^*$ continuous $\a$ but without no
continuity assumptions on $h$.
\begin{itemize}
\item [(i)]
The map $T: V \to C(X) , v \mapsto T(v)$, where $T(v): X \to \R$
is defined by $$T(v)(x)= \langle v, \alpha(x) \rangle$$  is a linear $S$-operator
(between right $S$-actions) with $||T|| \leq 1$.
\item [(ii)]
$T(v_0) \in \RUC_S(X)$ for every strongly  continuous
vector $v_0$ in $V$. Hence, if $h$ is strongly continuous
then $T(V) \subset \RUC_S(X)$.
%(respectively, $T(V) \subset WRUC_S(X)$ ).
\item [(iii)] If $V$ is reflexive, then $T(V) \subset \WAP_S(X)$.
\end{itemize}
\end{f}
\begin{proof} (i) is straightforward.

(ii) Observe that $||\a (x) || \leq 1$ for every $x \in X$. We get
$$||T(v_0)s - T(v_0)s_0|| = sup \{ |\langle v_0s -
v_0s_0, \a (x) \rangle|:  x \in X \}\leq$$
$$\leq ||v_0s - v_0s_0|| \cdot ||\a
(x)|| \leq ||v_0s - v_0s_0||. $$
This implies that $T(v_0) \in \RUC_S(X)$.

(iii) If $V$ is reflexive, the orbit $vS$ is relatively weakly
compact for each $v \in V$.
By the (weak) continuity of the $S$-operator $T$, the same
is true for the orbit of $T(v)$ in $C(X)$. Therefore we get $T(v) \in
\WAP_S(X)$.
\end{proof}

\begin{prop} \label{p-Telgen}
For every $S$-flow $X$ the following are equivalent:
\begin{enumerate}
\item
$f \in \RUC_S(X)$.
\item
There exist: a Banach space $V$, a strongly continuous
antihomomorphism $h: S \to \Theta(V)$, a weak$^*$ continuous
equivariant map $\a: X \to B^*$, and a vector $v \in V$ such that
$f(x)=\langle v, \a (x) \rangle$ (that is $f=T(v)$).
\end{enumerate}
\end{prop}
\begin{proof} (1) $\Longrightarrow$ (2) The function $f$
belongs to an $S$-invariant Banach subalgebra $\mathcal A$ of
$\RUC_S(X)$. The right action of $S$ on $V:=\mathcal A$ is jointly
continuous. Then by Fact \ref{Teleman}, corresponding left action
of $S$ on the dual ball $(B^*,w^*)$ is jointly continuous. Then
the naturally associated map $\a: X \to B^*$ and the vector $v:=f$
satisfy the desired property.

(1) $ \Longleftarrow$ (2) Immediate by Fact \ref{comesEb} (ii).
\end{proof}

%Actually it is easy to show that there exists a natural 1-1
%correspondence between jointly continuous $S$-compactifications of
%$X$, $S$-invariant Banach subalgebras $\mathcal A \subset
%RUC_S(X)$ and strongly continuous Banach representations of
%$(S,X)$. We postpone the details until Theorem \ref{} below where
%we deal also with reflexive and Asplund representations .

\begin{prop} \label{Tel2}
For every semitopological monoid $S$ the following are equivalent:
\begin{enumerate}
\item
$f \in \RUC(S)$.
\item
There exist: a Banach space $V$, a {\it strongly continuous}
antihomomorphism $h: S \to \Theta(V)$, and a pair of vectors $v \in
V$ and $\psi \in V^*$ such that $f=m_{v,\psi}$.
\end{enumerate}
\end{prop}
\begin{proof} (1) $\Longrightarrow$ (2)
Consider the Gelfand compactification $u_R: S \to S^R$ defined by
$\RUC(S)=C(S^R)$.  Then the action $S \times S^R \to S^R$ is
jointly continuous by Fact \ref{f-c=r}. Now define: $V:=C(S^R)$,
corresponding strongly continuous $h: S \to \Theta(V)$ (induced by
the right action of $S$ on $C(S^R)$), $v:=f \in V$ and $\psi=u_R(e)
\in V^*$.

(1) $ \Longleftarrow$ (2) Immediate by Fact \ref{introversion}.2.
\end{proof}

\vskip 0.1cm

As we already have seen a right uniformly continuous function can
be represented as a matrix coefficient $m_{v,\psi}$ of some
strongly continuous {\it Banach representation}. 
If $S$ is a group then $h(S) \subset \Is(V)$ in Proposition \ref{Tel2}. So, we can replace $\Theta(V)$ by $\Is(V)$.  

We mentioned 
the well known case of Hilbert representations. A positive
definite function on a topological group $G$ is exactly a matrix
coefficients of some unitary representation. One of our aims is to
understand the role of matrix coefficients for intermediate cases
of reflexive and Asplund representations.
%It turns out that the situation can be described adequately.
We show that wap functions
are exactly the {\it reflexive matrix coefficients}. In the
``Asplund case'' this approach leads to a definition of
{\it Asplund functions} introduced in
Section \ref{s-asp-f}.

\section{Reflexive representations of flows}

\vskip 0.1cm

\begin{defin} \label{duality}
A (bounded) {\it duality} is a separately continuous (resp., bounded)
mapping $\langle,\rangle: Y\times X \rightarrow \R$. We say that
the duality is {\it right strict} if the
corresponding continuous map 
$$q_X: X \to C_p(Y), \hskip 0.3cm
q_X(x)=\langle y,x \rangle$$ is a topological
embedding (e.g., an injection if $X$ is compact).

The ``left'' version can be defined analogously. Then ``strict''
will mean left and right strict simultaneously.

 Let a semigroup $S$ act
on $X$ and $Y$ by the following actions:
$$\pi_X: S \times X  \to X,
\hskip 0.3cm \pi_Y: Y \times S \to Y.$$
 The duality is an {\it $S$-duality} (or,
$S$-{\it invariant}) if
$
\langle ys,x\rangle =\langle y,sx\rangle 
$.
\end{defin}

Consider two typical examples:

\begin{enumerate}
\item ``Canonical reflexive duality'': $rB \times B^* \to \R$
with compact spaces $rB$ and $B^*$ (under weak topologies)
is defined for every reflexive $V$, a positive number $r>0$ and an antihomomorphism $h:
S\to \Theta(V)$. In particular, we can choose the natural action of
$S=\Theta(V)^{opp}$. Observe that $\Theta(V)^{opp}=\Theta(V^*)$
(compare Proposition \ref{ii}).
%mich
%\item ``Canonical reflexive duality'': $rB \times B^* \to [-1,1]$
%with compact spaces $B$ and $B^*$ (under weak topologies)
%is defined for every reflexive $V$ and an antihomomorphism $h:
%S\to \Theta(V)$. In particular, we can choose the natural action of
%$S=\Theta(V)^{opp}$. Observe that $\Theta(V)^{opp}=\Theta(V^*)$
%(compare Proposition \ref{ii}).
\item Let $K\subset V$ be a weakly compact S-invariant subset in a
Banach space $V$ with respect to some antihomomorphism $h: S \to
\Theta (V)$. Then $K \times B(V^*) \to \R$ is a left strict
$S$-duality.
\end{enumerate}

\begin{lemma} \label{net}
Let $\langle ,\rangle : Y \times X \to \R$ be an $S$-duality.
\ben
\item
$\langle ,\rangle $ is left strict iff a net
$y_i$ converges to $y$ in $Y$ exactly when
$\langle y_i,x\rangle  \to \langle y,x\rangle $ in $\R$ for every $x \in X$. Similarly,
$\langle ,\rangle $ is right strict iff a net $x_i$ converges to $x$ in $X$
exactly when $\langle y,x_i\rangle  \to \langle y,x\rangle $ in $\R$ for every $y \in Y$.
\item
 Let $\langle ,\rangle $ be a left (right)  strict $S$-duality.
Then all $s$-translations  ${\tilde s}: Y \to Y$
(resp., ${\tilde s}: X \to X$)
are continuous.
\item
 Let $\langle ,\rangle $ be a strict $S$-duality.
Then $\pi_Y$ is separately continuous iff $\pi_X$ is separately continuous.
\item
Let
%(separately continuous)
$\langle ,\rangle : Y \times X \rightarrow \R$ be a left strict $S$-duality.
Then it can be
reduced canonically to the naturally associated strict
%(separately continuous)
$S$-duality $\langle ,\rangle _q: Y \times X_q \rightarrow \R$.
%Moreover, $\langle ,>_q=\langle ,>$ iff $\langle ,>$ is separated.
If $\pi_Y$ is separately continuous then $\pi_{X_q}$
is also separately continuous.
\een
\end{lemma}
\begin{proof} (1) Follows from the net characterization of
the product topology.

(2) We have to show that every $s$-translation $\pi_Y^s: Y \to Y$
is continuous (the case of $\pi_X^s$ is similar). Let
$y_i \to y$. In order to show that
%mm22Jan07 here we need   ys  and not sy
$sy_i \to sy$, it suffices by (1) to check that $\langle y_i s, x\rangle   \to
\langle ys,x\rangle $ for each $x \in X$. Or, equivalently, we have to show that
$\langle y_i,sx\rangle   \to \langle y,sx\rangle $. The latter follows from the assumption
$y_i \to y$ and the separate continuity of $\langle ,\rangle $.

(3) Similar to the proof of (2).

(4) Consider the canonical continuous map
$q_X: X \to C_p(Y)$ and the corresponding range $X_q=q_X(X)
\subset C_p(Y)$. Define
$$\langle ,\rangle _q: Y \times X_q \rightarrow \R, \hskip 0.1cm \langle y, q_X(x)\rangle _q := \langle y,x\rangle .$$
It is easy to show that this is a well-defined strict duality.
\delete{ It is more convenient verify this assertion in two steps.
  First "reducing from the right"  consider the following duality

\centerline{$ Y \times X_q \rightarrow \R, \hskip 0.1cm \langle y,
q_X(x)\rangle  :=\langle y,x\rangle $}

Clearly this duality is right strict and $S$-invariant. Now, as a
second step, "reducing from the left side" similarly we get the
duality \centerline{$ {\tilde Y} \times X_q \rightarrow \R, \hskip
0.1cm \langle {\tilde q_Y}(y), q_X(x)\rangle  :=\langle q(y),x\rangle $}

Here ${\tilde Y}$ is the subspace ${\tilde q_Y}(Y)$ of $C_p(X_q)$
and ${\tilde q_Y}: Y \to C_p(X_q)$ is the natural continuous
mapping induced by the ``first step duality''. Consider also the
{\it injective} map $\nu: C_p(X_q) \to C_p(X)$ induced by the
surjective map $q:_X: X \to X_q$. Injectivity of $\nu$ implies
that ${\tilde q_Y}(y)=q_Y(y)$. We can identify ${\tilde q_Y}(X_q)$
with $\nu ({\tilde q_Y}(Y))=q_Y(Y)$ not only as a set but even as
a topological spaces . Indeed , observe that for every net $y_j$
in $Y$ we have the following equivalence:
 $\langle y_j,x\rangle  \to \langle y,x\rangle $ for every $x \in X$ iff
$\langle y_j,q(x)\rangle  \to \langle y,q(x)\rangle $ for every $x \in X$.

Thus, $\langle ,\rangle _q$ is a well-defined duality. By the construction it is
a strict duality} Moreover, the action of $S$ on $X$ induces the
natural action of $S$ on $X_q$ such that $q: X \to X_q$ is
$S$-equivariant and $\langle ,\rangle _q$ is $S$-invariant. Apply (3) to
$\langle ,\rangle _q$. If $\pi_Y$ is separately continuous then $\pi_{X_q}$ is
also separately continuous by virtue of (3). 
%This is trivial for $\pi_{X_q}$ because
%$X \to X_q$ is continuous. Now continuity of orbit maps for
%$\pi_{Y_q}$ follows from (3).
% and $Y \to Y_q$ are continuous $S$-maps.
\end{proof}

\begin{prop} \label{ii}
Let $X$ and $Y$ be Hausdorff $S$-flows such that $X$ is
compact and every orbit closure $cl(yS)$ in $Y$ is compact. Assume that
$\langle ,\rangle : Y\times X \to \R$ is
a strict $S$-duality. Then the Ellis semigroup $E(S,X)=E(X)$
is semitopological and there exist separately continuous actions
$E(X) \times X \to X$ and $Y \times E(X) \to Y$ extending
the original actions of $S$ and
such that $\langle ,\rangle : Y\times X \to \R$ becomes an $E(X)$-duality .

%with the canonical semigroup
%compactifications $\lambda: S \to E(X)$ and $\varsigma: S \to
%E(Y)$.
%\ben
%\item
%There exists a topological isomorphism $\gamma: E(X) \to E(Y)$
%such that $\lambda \circ \gamma = \nu$ and
%$<,>$ becomes an $E(X)$-duality with respect to
%the natural action $E(X) \times X \to X, \hssk px=p(x)$ and the
%right action $Y \times E(X) \to Y, \hssk yp=\gamma(p)(y)$.
%%Moreover, $\lambda(s)x=sx$ and $y \lambda(s)=y \varsigma(s) = ys$
%%for every $(s,x,y) \in S \times X \times Y $.

%\item
% $E(X)$ and $E(Y)$ are compact semitopological semigroups and
%the pairs $(E(X),X)$ and $(Y,E(Y))$ are semitopological wap flows.
%\een
\end{prop}
\begin{proof}
%%may29
%We have to show that $p: X \to X$ is continuous for every $p \in
%E(X)$.
By Lemma \ref{net} (assertions (2) and (3)) we need only to check
that there exists an action of $E(X)$ on $Y$ which extends the
original right action of $S$ on $Y$ in such a way that the duality
$\langle ,\rangle $ becomes an $E(X)$-invariant.

Denote by $\lambda: S \to E(S,X)$ the canonical semigroup
compactification. We follow the idea of \cite[Proposition
II.2]{EN}. Let $p \in E(X)$. Choose arbitrarily a net $s_i \in S$
such that $\lambda (s_i)$ converges to $p$ in $E(X)$. Let $y \in
Y$. Using the compactness of $cl(yS)$, one can pick a subnet $t_j$
of $s_i$ such that
%%may29
%%$\varsigma (t_j)(y)=
$y t_j$ converges to some $z\in Y$. Define $yp:=z$. Then for every
$x \in X$ we have
$$\langle z,x\rangle =\langle \lim(y t_j), x\rangle =\lim \langle y t_j, x\rangle =\lim \langle y, t_j x\rangle =$$
$$=\langle y, \lim (t_j x)\rangle = \langle y, \lim \lambda(t_j)x\rangle =\langle y, px\rangle .$$
 The element $\langle y, px\rangle $ does not depend on the choice of
subnets in the definition of $z$.
Since the duality $Y \times X \to \R$ is (left) strict, we can
conclude that such $z$ is uniquely determined. Thus,
%%may29
%$z$ and hence also
$yp$ is well-defined. These computations show also that $\langle yp,
x\rangle =\langle y, px\rangle $ because each of them is the limit of $\langle yt_j, x
\rangle =\langle y, t_j x\rangle $.
%%may29
%%Clearly, this action of
%$E(Y)$
Since the given duality is strict it follows that the function $Y
\times E(X) \to Y, \quad (y,p) \mapsto yp$ is a right action which
extends the given action of $S$ on $Y$. Indeed, we can choose for
$p:=\lambda(s)$ the constant net $s_i=s$ in the above definition.
\end{proof}

\begin{remarks} \label{r-two} \ 
\bit \item [(a)] If $Y$ is compact then it is easy to see that the
Ellis semigroups $E(S,X)$ and $E(Y,S)$ are antiisomorphic as
semitopological semigroups.

\item [(b)] Proposition \ref{ii} provides a proof of the crucial
implication (1) $\Longrightarrow$ (2) of Fact \ref{EN}. Indeed,
endow the $S$-flow $Y:=C(X)=\WAP_S(X)$ with the pointwise topology
and apply Proposition \ref{ii} to the natural $S$-pair $Y \times X
\to \R$.

\item [(c)]
As a corollary we get that the semitopological flow
$(E(X),X)$ of Proposition \ref{ii} is wap. In fact this flow is
even Eberlein as it follows by Theorem \ref{Davis} below.

\item [(d)]
Proposition \ref{ii} implies that a compact $S$-flow $X$ is
$S$-wap iff $B^*=B(C(X)^*)$ is $S$-wap. Indeed, by Fact \ref{EN}
it suffices to show that the action of the Ellis semigroup
$E(S,X)$ on $B^*$ is separately continuous. This follows from
Lemma \ref{net}.3 because by Fact \ref{f-functions} (i) the action
of $E(S,X)$ on $B(C(X))_w$ is separately continuous and $B(C(X))_w
\times B^*_{w^*} \to \R$ is a strict $E(S,X)$-duality. Therefore
if $X$ is a compact $S$-wap then we can conclude that the
$S$-subflow $P(X)\subset B^*$ of all probability measures is wap,
too. This fact was established earlier by Glasner \cite{Gl}. Below
we show (Theorem \ref{AL}) that $P(X)$ is $S$-Eberlein if $X$ is
$S$-Eberlein. \eit
\end{remarks}

\vskip 0.1cm Now we prove that all bounded strict $S$-dualities $Y\times X \to
\R$, with compact $X,Y$, come as restrictions of canonical reflexive $S$-dualities $rB \times B^* \to \R$.  
This fact is interesting even in the purely topological
context (that is, for trivial $S$) which has been proved by
Krivine and Maurey \cite{KM} for metrizable compacta $X$ and $Y$. 
In the proof we provide a modification for flows of the well-known
construction of Davis, Figiel, Johnson and Pelczynski \cite{DFJP}, taking into account 
results of Lima, Nygaard and Oja \cite{LNO}, getting an operator $j: V \to E$ (see proof of Theorem \ref{Davis}) which is not only bounded but even $||j|| \leq 1$. This allows us to get a little bit more than in \cite{KM}. 
Namely, representing every strict duality $Y\times X \to
[-1,1]$ as a restriction of the canonical reflexive duality $B \times B^* \to
[-1,1]$ (with $r=1$). 

%  
%\vskip 2cm 
%%160919
%\begin{remark} \label{r:iso-DFJP}
%	The fundamental DFJP-factorization construction from \cite{DFJP} has an ``isometric modification".
%	According to \cite{LNO} one may assume in Theorem \ref{t:general}
%	that the bounded operator $j: V \to \mathcal{A}$ has the property $||j|| \leq 1$.
%	More precisely, we can replace in the Equation \ref{bLNO}
%	the sequence of sets $M_n:=2^n W + 2^{-n} B$ by $K_n:=a^{\frac{n}{2}} W + a^{-\frac{n}{2}} B$, where
%	$2 < a < 3$ is the unique solution of the equation $\sum_{n=1}^{\infty} \frac{a^n}{(a^n+1)^2} =1$.
%	For details see \cite{LNO}.
%	Taking into account this modification (which is completely compatible with our $S$-space setting)
%	for a set $F \subset C(X)$ with $\sup \{|f(x)| : \ x \in X, f \in F\} \leq 1$
%	we can assume that $\nu(F) \subset B$ and $\a(X) \subset B^*$. Hence the following sharper diagram commutes
%	
%	\begin{equation} \label{diag2}
%	\xymatrix{ F \ar@<-2ex>[d]_{\nu} \times X
%		\ar@<2ex>[d]^{\a} \ar[r]  & [-1,1] \ar[d]^{id} \\
%		B \times B^* \ar[r]  &  [-1,1] }
%	\end{equation}
%	Note also that this modified version from \cite{LNO} of the DFJP-construction repairs in particular
%	the proof of \cite[Theorem 4.5]{Me-nz}.
%	% The latter was first corrected in the arxiv version of \cite[Theorem 4.5]{Me-nz}
%	%using, however, diagrams like \ref{diag1}, where $\nu(F)$ and $\a(X)$ are bounded.
%\end{remark}
%%%%%%%%% !!!!!!!!!!!!!!!!!
%\vskip 2cm 

\begin{thm} \label{Davis}
Let $\nu: Y \times X \rightarrow [-r,r]$ be a bounded strict
$S$-duality, with semitopological compact $S$-spaces $X$ and $Y$ for some positive number
$r >0$.
Then there exist: a reflexive Banach space $V$, a weakly continuous antihomomorphism $h: S \to \Theta(V)$ 
and weak
%%may29  bounded
embeddings $\gamma_1: X \to B^*$ and $\gamma_2: Y \to B$ such
that the following diagram is commutative.
%nd
%\begin{diagram}[height=2.2em,width=3.0em,tight]
% &Y \times X &\rTo^{} & \R  \\
% &\dTo^{\gamma_2} \quad \quad \dTo_{\gamma_1}  &    &\dTo_{id}  \\
% &rB \times B^* & \rTo^{} & \R \\
%%mich
%%Every strict $S$-duality $Y \times X \rightarrow [-1,1]$, with
%%semitopological compact $S$-spaces $X$ and $Y$,
%%is an $S$-restriction of a reflexive
%%$S$-duality $B \times B^* \rightarrow [-1,1]$ which is induced by a
%%weakly continuous antihomomorphism $h: S \to \Theta(V)$.
%%%(see the diagram).
%%\begin{diagram}[height=2.2em,width=3.0em,tight]
%% &Y \times X &\rTo^{} & [-1,1]  \\
%% &\dTo^{\gamma_2} \quad \quad \dTo_{\gamma_1}  &    &\dTo_{id}  \\
%% &B \times B^* & \rTo^{} & [-1,1] \\
%\end{diagram}

\begin{equation*} \label{diag2}
\xymatrix{
Y \ar@<-2ex>[d]_{\gamma_2} \times X
\ar@<2ex>[d]^{\gamma_1} \ar[r]  & [-r,r] \ar[d]^{id} \\
rB \times B^* \ar[r]  &  [-r,r]}
\end{equation*}

If $S=G$ is a group then $h(G) \subset \Is(V)$. 
If the action $S \times X \to X$ is jointly continuous then we can
suppose that $h: S \to \Theta(V)$ is strongly continuous.
 \end{thm}
\begin{proof} 
	It is enough to prove for the particular case of $r=1$. 
By Proposition \ref{ii} we can suppose that $S$ is a compact
semitopological semigroup. Adjoining the isolated identity $e$,
one can assume
even that $S$ is a monoid and $e$ acts as the identity map. The map $q_Y : Y
\to (C(X), w)$ is a topological embedding by Grothendieck's Lemma.
We will identify $q_Y(Y)$ and $Y$.
%(use: $\gamma(Y)$ is bounded
%and $p$-compact).
Denote by $E$ the Banach subspace of $C(X)$ topologically
generated by $Y$. That is, $E=cl(sp(Y))$.
%We will identify $Y$ and $\gamma_2 (Y)$.

Consider the right action $C(X) \times S \rightarrow C(X), \hskip
0.1cm (fs)(x):=f(sx)$. Then every $s$-translation ${\tilde s}:
C(X) \to C(X)$ is a contractive linear operator. The orbit map
${\tilde y}: S \to C(X)$ is $p$-continuous for every $y \in Y
\subset C(X)$. By our assumption $S$ is compact. Therefore, the
orbit $yS$ is bounded $p$-compact, and, hence $w$-compact by
Grothendieck's Lemma. Since $p$-topology coincides with the
$w$-topology on $yS$, it follows that ${\tilde y}: S \to C(X)$ is
$w$-continuous. Then the same is true for every $u \in
E=cl(sp(Y))$ (as it follows, for example, from \cite[Proposition
6.1.2]{BJM}). By the
 Hahn-Banach Theorem, the weak topology of $E$ is the same as its
relative weak topology as a subset of $C(X)$. Therefore we get
that $((E,w), S)$ is a semitopological flow. Consider the convex
hull $co(-Y\cup Y)=W$. By the Krein-Smulian Theorem, $W$ is
relatively weakly compact (convex and symmetric) subset in $E$.
%mich    %Since $A$
Since $r=1$, $||q_Y(y)|| \leq 1$ (in $E$) for every $y \in Y$. So, in addition, we get that $W$ is a subset of the unit ball $B(E)$. 

We apply an ``isometric modification" \cite{LNO} of the fundamental DFJP-factorization construction from \cite{DFJP}.  
For each natural $n$, set $U_n = a^{\frac{n}{2}} W + a^{-\frac{n}{2}} B(E)$, where
$2 < a < 3$ is the unique solution of the equation $\sum_{n=1}^{\infty} \frac{a^n}{(a^n+1)^2} =1$ (for details see \cite{LNO}).

 Let $\| \
\|_n$ be the Minkowski's
functional of the set $U_n$. That is, $\| v\|_n = \inf\
\{\lambda > 0 \bigm| v\in \lambda U_n\}$.

\nt Then $\| \ \|_n$ is a norm on $E$ equivalent to the given norm
of $E$. For $v\in E,$ let
$$N(v):=\left(\sum^\infty_{n=1} \| v \|^2_n\right)^{1/2} \hskip 0.1cm \text{and}
\hskip 0.1cm \hskip 0.1cm V: = \{ v \in E \bigm| N(v) < \infty
\}.$$

 Denote by $j: V \to E$ the inclusion map. Then:

(1) $(V, N)$ is a reflexive Banach space, $j: V \to E$ is a
%mich
%continuous linear injection, $Y \subset W \subset B(V)$ and
%$||j|| \leq 1$.
continuous linear injection, $Y \subset W \subset B(V)$ 
and $||j|| \leq 1$. 

%Indeed, first observe that $W \subset B(E)$.
%This implies that if $z \in Z$ and $||j(z)||=1$
%then $||z||_n$ is at least $2^{-n}$ for every natural $n$. Therefore
%$N(z)$ is at least 1.
\vskip 0.1cm

(2) The restriction of $j: V \ri E$ on each bounded subset $A$ of
$V$ induces a homeomorphism of $A$ and $j(A)$ in the weak
topologies.

%\delete{
%\proof \ \
%Consider the weak closure $cl_w(A)$ of $A$ in $X_f.$
%By the reflexivity of $X_f,$ the set $cl_w(A)$ is weakly compact.
%Hence, $j$, being weakly continuous and injective,
%induces a homeomorphism of
%$cl_w(A)$ and $j(cl_w(A))$ with respect to the weak topologies.
%}

\vskip 0.1cm

By our construction $W$ and $B(E)$ are $S$-invariant. Thus we get

(3)  $V$ is an $S$-subset of $E$ and $N (vs) \leq N(v)$ for every
$v\in V$ and every $s\in S$.

%\delete{
%\proof \ \
%It suffices to show that $\| vs \|_n \leq \| x \|_n$ for every $n\in
%\N.$
%By our construction $W$ and $B_E$ are $S$-invariant, that is
%$Ws \subset W$ and $B_E s \subset B_E$.
%Then, from $v\in \lambda(2^n W + 2^{-n} B_E))$
%we obtain that $vs \in \lambda(2^n (Ws) + 2^{-n} (B_E)s \subset
%\lambda(2^n W + 2^{-n} B_E).$
%Hence, $\|v s \|_n \le \| v\|_n,$ as required.
%%This proves assertion (6).
%}

\vskip 0.1cm

%As a corollary, we get that $Z$ is an $S$-invariant subset of $E$.
%Therefore, the restricted action
%$\a_f : S\times X_f \ri X_f$
%is well-defined.

(4) For every $v\in V$, the orbit map $\tilde{v}  : S\ri V, \hskip
0.2cm \tilde{v}  (s) = vs$ is weakly continuous.

\vskip 0.1cm Indeed, by (3), the orbit $\tilde{v}  (S) = vS$ is
$N$-bounded in $V$. Our assertion follows from (2) (for $A=vS$),
taking into account that $\tilde{v} : S \to E$ is weakly
continuous.

\vskip 0.1cm

By (3), for every $s\in S$, the translation map ${\tilde s}: V \to
V, v \mapsto vs$ is a linear contraction of $(V, N)$. Therefore,
we get the antihomomorphism $h : S \to \Theta(V), \hskip 0.3cm h
(s) = {\tilde s}$.

\vskip 0.1cm

Now, directly from (4), we obtain the following assertion.

(5)\ \ $h : S \to \Theta(V)$ is a $w$-continuous monoid
antihomomorphism.

\vskip 0.1cm By (1) and (2) the natural inclusion map $\gamma_2: Y
\to B=B(V)$ is a topological (weak) embedding. Define the
weak star embedding $\gamma_1: X \to V^*$ by
$\gamma_1(x)(v)=j(v)(x)=\langle v,x\rangle $. Clearly, $\gamma_1(x)=(i^* \circ j^* )(x)$,
where $i^*: C(X)^* \to E^*$ is the adjoint of
the inclusion $i: E \to C(X)$. Clearly, $||i|| = 1$. 
Since $||j|| \leq 1$ (by (1)) we get that
$\gamma_1(X)$ is a subset of the unit ball $B^*$ of $V^*$.

It is evident that $\gamma_1$ is
$S$-equivariant and weak$^*$ (=weak) continuous. On the other
hand, we have $\langle \gamma_2(y),~\gamma_1(x)\rangle =\langle y,x\rangle $. Since the original
duality is strict, we obtain that $\gamma_1$ is injective and
hence a topological (weak) embedding. 

\vskip 0.1cm  If $S \times X \to X$ is jointly continuous, then
the action $C(X) \times S \to C(X)$ is jointy continuous with
respect to the norm. By the definition of the Banach space
$(V,N)$, it is straightforward to show that all orbit maps
${\tilde v}: S \to V$ are $N$-norm continuous (recall that each
$\| \ \|_n$ is equivalent to the norm of $E$). This will guarantee
that $h: S \to \Theta(V)_s$ is continuous.
\end{proof}

Now we can prove the representation theorem.

\begin{thm} \label{main}
({\bf WAP Representation Theorem}) Let $X$ be a semitopological
$S$-flow. The following conditions are equivalent:
\begin{itemize}
\item [(i)]  $f : X \to \R$ is weakly almost periodic.
\item [(ii)] There exist: a representation $(h,\alpha)$ of $(S,X)$ into
reflexive $V$ with a weak continuous antihomomorphism $h: S \to
\Theta(V)$, weak continuous $\alpha: X \to B^*$, and a vector $v
\in V$ such that $f(x) =\langle v, \alpha(x)\rangle $.
\item [(iii)] As in (ii) but with no continuity assumptions on $h$.
\end{itemize}

If either: a) $S=G$ is a semitopological group; or b) $X$ is
compact and the action $S \times X \to X$ is jointly continuous,
then in (ii) we can suppose that $h$ is strongly continuous.
\end{thm}
\begin{proof} (i) $\Longrightarrow $ (ii)
We can suppose that $S$ is a monoid.
%mich
%and $\|f\| \leq 1$.
For the
desired representation of $f\in \WAP_S(X)$ by some reflexive $V$,
choose a left strict bounded duality
$$
\langle ,\rangle : K_f \times D_f \rightarrow \R
$$
where $K_f=(cl_w(fS), w)$ and  $D_f=(B(\WAP_S(X)^*),w^*)$.
The weak and pointwise topologies coincide on $K_f$. Therefore
the action of $S$ on $K_f$ is separately continuous. Note, however,
that the action of $S$ on $D_f$ is not necessarily separately
continuous. By Lemma \ref{net} we can pass to the naturally
associated strict separately continuous $S$-duality $\langle ,\rangle _q: K_f
\times (D_f)_q \rightarrow \R$.
%The action $K_f \times S \to K_f$ is separately continuous.
Lemma \ref{net}.3 guarantees
that the action of $S$ on $(D_f)_q$ is also separately continuous.
Denote by $t: X \to (D_f)_q$ the composition of two natural maps
$X \to D_f$ and $D_f \to (D_f)_q$. Now, by Theorem \ref{Davis},
there exist: a reflexive Banach space $V$ and a weakly continuous
antihomomorphism $h: S \to \Theta(V)$ such that $S$-duality
$\langle ,\rangle_q$ equivariantly can be realized as a part of a reflexive
%mich
duality $rB \times B^* \rightarrow \R$.
%where $\|\phi\|\leq n$ for every $\phi \in K_f$. We get the diagram.
%nd
%\begin{diagram}[height=2.2em,width=3.0em,tight]
% & K_f \times (D_f)_q &\rTo^{} & \R  \\
% &\dTo^{\gamma_2} \quad \quad \dTo_{\gamma_1}    &  &\dTo_{id}  \\
% &rB \times B^* & \rTo^{} & \R \\
% \end{diagram}

\begin{equation*} \label{diag3}
\xymatrix{
K_f \ar@<-2ex>[d]_{\gamma_2} \times (D_f)_q
\ar@<2ex>[d]^{\gamma_1} \ar[r]  & \R \ar[d]^{id} \\
rB \times B^* \ar[r]  &  \R }
\end{equation*}

Define $\alpha: X \to B^*$ by
$\alpha(x)=\gamma_1(t(x))$ and
pick $v:=\gamma_2(f)$. Then $f(x) =\langle v,\alpha(x)\rangle $, as desired.

 (ii) $\Longrightarrow $ (iii) Is trivial.

(iii) $\Longrightarrow $ (i) Is immediate by Fact \ref{comesEb}
(iii).

If $S$ is a semitopological group, then every weakly continuous
reflexive (anti)representation is automatically strongly
continuous as we mentioned in Section \ref{s-repr} (see Remark
\ref{r-automat}). This proves the case ``a)''. In the second case
``b)'', we can apply directly Theorem \ref{Davis}.
\end{proof}

Now we easily obtain one of our main results.

\begin{thm} \label{wapeb}
An $S$-flow $X$ is wap iff $X$ is REFL-approximable.
\end{thm}
\begin{proof} If $X$ is REFL-approximable then $X$
is wap  by Fact \ref{comesEb} (iii).

The nontrivial part follows from Theorem \ref{main}
because if $X$ has sufficiently many wap functions, then $(S,X)$
has sufficiently many reflexive representations.
\end{proof}

\begin{corol} \label{c-wap-RN}
Every wap flow $(S,X)$ is RN-approximable.
\end{corol}

It is well known that a countable product of Eberlein (RN) compacta
is again Eberlein (resp.: RN). We show that the same is true for flows.

\begin{lemma} \label{countpr}
The classes of Eberlein and RN $S$-flows are closed under countable products.
\end{lemma}
\begin{proof} Let $X_n$ be a sequence of Eberlein (or, RN)
$S$-flows.  By the definition there exists a sequence of reflexive
(Asplund) representations
$$ (h_n,\a_n): (S,X_n) \rightrightarrows
(\Theta(V_n)^{opp},B(V_n^*)).$$
We can suppose that each $X_n$ is compact
and
%corrected August 21 2005
%$h(X_n) \subset 2^{-n}B(V_n)$. Turn to the $l_2$-sum of
$\a_n(X_n) \subset 2^{-n}B(V_n^*)$. Turn to the $l_2$-sum of
representations. That is, consider
$$ (h,\a): (S,X) \rightrightarrows (\Theta(V)^{opp}, B(V^*))$$
where $V:=(\sum_n V_n)_{l_2}$, $h(s)(v)=\sum_n h(s)(v_n)$ for
every $v=\sum_n v_n$, and $\a(x)= \sum_n \a_n(x_n)$ for every
$x=(x_1, x_2, \cdots) \in \prod_n X_n$.  It is easy to show that
$\a(x) \in B(V^*)$, $\a$ is weak$^*$ continuous and injective
(hence, a topological embedding). Now use the fact that the
$l_2$-sum of reflexive (Asplund) spaces is again reflexive
(Asplund) \cite{Fa}.
\end{proof}

\begin{corol} \label{seccwap} \ 
\begin{itemize}
\item
[(i)] Every second countable wap flow is Eberlein.
\item
[(ii)] Every second countable RN-approximable flow is an RN flow.
\end{itemize}
\end{corol}
\begin{proof}
Assertion (ii) is immediate by Lemma \ref{countpr}.
For (i), we need also Theorem \ref{wapeb}.
\end{proof}

%\delete{
%\begin{prop}
%The following are equivalent:
%
%(i) A flow (S,Y) is Eberlein
%(there exists a faithful reflexive representation $Y \to B_V$)
%
%
%(ii) There exists a proper $S$-compactification $Y \to D$ and a
%right strict bounded $S$-duality $K \times D \to \R$.
%
%(iii) There exists a countable
%right strict family $w_n: X_i \times Y \ri \R$
%$(n \in \N)$ of (left strict) $S$-dualities.
%
%(iv) There exists a countable
%topologically separating family $m_n: (S,X) \ri (\Theta(V_n), B_{V_n})$
%$(n \in \N)$ of reflexive representations.
%
%(v) $X$ is a relatively weakly compact subflow in a Banach $S$-space.
%
%\end{prop}
%}

%\begin{proof}
%In the part $(iv) \Longrightarrow (v)$
%Use $l_2$-sum and $2^{-n}B_{V_n}$ balls.
%\end{proof}

The following theorem provides, in particular, a
flow generalization of a result by
Amir-Lindenstrauss \cite{AL}
which states that if $X$ is an Eberlein compact then
$B^*=(B(C(X)^*), w^*)$ is Eberlein, too.

\begin{thm} \label{AL}
Let $X$ be a compact semitopological $S$-flow.
The following are equivalent:
\begin{itemize}
\item
[(i)] $X$ is $S$-Eberlein.
\item
[(ii)] There exists a Banach space $E$, a homomorphism $h: S \to \Theta(V)$
(no continuity assumptions on $h$), and an $S$-embedding $\a: X \to (V,w)$.
\item
[(iii)] There exists a compact space $Y$ and a (right) strict
$S$-duality $Y \times X \to \R$ .
\item
[(iv)] There exists a sequence of $S$-invariant weakly
compact subsets $K_n \subset C(X)$ such that $\cup_{n \in \N} K_n$
separates the points of $X$.
\item
[(v)] There exists an $S$-invariant weakly compact subset $M$ in
$C(X)$ such that $cl(sp(M))=C(X)$.
\item [(vi)] $B^*$ is $S$-Eberlein.
\item [(vii)]
 $P(X)$ is $S$-Eberlein.
\end{itemize}
\end{thm}
\begin{proof}
(i) $\Longrightarrow$ (ii) By the definition there exists a faithful
reflexive $V$-representation. That is, we can choose a weakly continuous
homomorphism $h: S \to \Theta(V)^{opp}=\Theta(V^*)$ and an equivariant
embedding $\alpha: X \to B(V^*)$. It suffices to choose $E:=V^*$.

(ii) $\Longrightarrow$ (iii)
By our assumption $X$ is $S$-embedded into $(E,w)$. Define the
right strict $S$-duality $Y \times X \to \R$
as a restriction of the canonical duality
where $Y:=(B(E^*), w^*)$.

(iii) $\Longrightarrow$ (iv)
Use the ``right version'' of Lemma \ref{net}.4. Then
our right strict $Y \times X \to \R$ duality
induces the strict $S$-duality $Y_q \times X \to \R$.
We can suppose in addition that this
duality is bounded. Now define simply
$K_n:=Y_q \subset C(X)$ for each $n$ and use Fact \ref{Groth}.

(iv) $\Longrightarrow$ (v) Look at $K_n$ as an $S$-subflow of
$(C(X),w)$. We can suppose that $K_n\subset B(C(X))$. Following
a method of Rosenthal \cite{Ro1}, consider $S$-invariant
set $M_n$ consisting
of the constant function equal to 1 on $X$ and of all products of
functions $f_1 \cdot f_2 \cdots f_n$ where $f_i \in (\cup_{m=1}^n
K_m) \cup \{1\}$. By Fact \ref{Groth} and the Eberlein-Smulian
theorem, it is easy to see that each $M_n$ is weakly compact.
Then $M: =\cup_{n\in \N} 2^{-n} M_n$ is also
$S$-invariant and weakly compact in $C(X)$.
By the Stone-Weierstrass theorem, $sp(M)$ is dense in $C(X)$.

(v) $\Longrightarrow$ (vi) We can suppose that
$M \subset B(C(X))$. The corresponding left strict
$S$-duality $M \times B^* \to [-1,1]$ is also right strict because
$cl(sp(M))=C(X)$. Now
we can apply Theorem \ref{Davis}.

%As in the proof of Theorem \ref{main} one can construct a
%reflexive Banach space $V$ with a separately w-continuous linear
%right action of $S$ and a natural continuous linear $S$-injection
%$j: V \to C(X)$. By assertion (1), this map has a dense range.
%Therefore the adjoint $j^*: C(X)^* \to V^*$ is an injective
%$S$-map. This adjoint is weak$^*$-weak continuous (reflexivity of
%$V^*$). Therefore
% $B(C(X)^*)$ canonically is $S$-homeomorphic to a
%weak compact subset of $V^*$. Thus we get a faithful
%representation of the flow $(S, B^*)$ in reflexive $V$.

(vi) $\Longrightarrow$ (vii) and (vii) $\Longrightarrow$ (i)  are
trivial because $P(X)$ is an $S$-subflow of of $B^*$ and $X$ can be treated
as an $S$-subflow of $P(X)$. 
%\delete{
%In particular, $sp(A)$ separates points of $C(X^*)$. Therefore,
%even $A$ separates the points of $C(X^*)$. In particular, $A$ separates
%points of $B^*$. By Lemma \ref{distel} every
%$Y^{(n)}$ is $w$-compact . Moreover $Y^{(n)}$ is $S$-flow and
%the natural pairing $w_n: Y^{(n)} \times B^* \to \R$ is a
%left strict $S$-duality.
%Since the corresponding sequence of reflexive
%representations $\alpha_n: B^* \to B(V)_n$
%separates points of $B^*$ (we can omit ${\bf 1}$),
%then we obtain that $B^*$ is $S$-wap.
%
%(ii)
%By Theorem \ref{dualities} there exists
%a $X$-separating
%family $\Sigma=\{Y_i \times X\to \R | i \in I \}$ of left strict
%$S$-dualities. Then
%$\{B^* \times (K_{i_1} \cdot K_{i_2} \cdots K_{i_n}) \to \R |
%\{i_1, i_2, \cdots, i_n\} \subset I\}$ is a $B^*$-separating family of
%right strict $S$-dualities.
%}
\end{proof}

%\delete{
%\begin{corol} \label{measures}
%$X$ is $S$-Eberlein ($S$-wap) iff $M(X)$ (space of probability measures)
%is $S$-Eberlein ($S$-wap).
%\end{corol}
%
%The fact $X \in WAP$ iff $P(X) \in WAP$ was established earlier
%by Glasner \cite{Gl}.
%}

\section{Reflexive representations of (semi)groups}

\vskip 0.1cm

 Now we examine a particular but important case of the flows
$(S,S)$, left regular actions of a semitopological semigroup
$S$ on itself by multiplication. Every compact semitopological
semigroup is wap. In  general, $S$ is wap iff the universal
semitopological compactification $S \to S^W$ is an embedding iff
$S$ is a subsemigroup of a compact semitopological semigroup.

Every locally compact Hausdorff topological group $G$ is wap being
a subsemigroup of its one-point compactification (which clearly is
a compact semitopological semigroup). Moreover, it is well known that
such $G$ is even unitarily representable
because it can be embedded into the unitary group $\Is(H)_s$
of the Hilbert space $H=L_2(G, m_{Haar})$, where $m_{Haar}$ is
the Haar measure on $G$.

It is also easy to show that every {\it non-Archimedean}
(having a local base of open subgroups) topological group is
unitarily representable (and, hence wap).
Distinguishing unitarily and reflexive representability
(and answering a question of Shtern \cite{Sh}),
we show in \cite{Meist} that the additive group of
$L_4[0,1]$ is wap but not unitarily representable.
The proof is based on Grothendieck's double limit property for wap
functions.
It is still an open question if every abelian Hausdorff topological group
(e.g., the additive group of a
Banach space) is wap.

Not every topological (even Polish) group is wap. Indeed, the
group $G=\Homeo_+[0,1]$ of all orientation preserving
selfhomeomorphisms of the closed interval is not wap \cite{Merup}.
In fact we show that every wap function on such $G$ is necessarily
constant (conjectured by Pestov). As a corollary this implies that
the universal semitopological compactification $G^W$ of $G$ is
trivial (answering a question of Ruppert \cite{Ru}) and every
(weakly) continuous bounded representation $h: G \to Aut(V)$ into a
reflexive space $V$ is trivial. This example also shows (answering
a question of Milnes \cite{Mi}) that there exists a nonprecompact
Hausdorff topological group $G$ such that $\WAP(G)=\AP(G)$.

\vskip 0.1cm Turn again to the WAP Representation theorem. It
implies that every wap function comes from a reflexive matrix
coefficient.

\begin{thm} \label{matrixwap}
For every semitopological
monoid $S$ the function $f: S \to \R$ is wap iff $f$ is a
matrix coefficient of a weak continuous
antihomomorphism $S \to \Theta(V)$ for a reflexive $V$.
That is, there exist $v \in V$ and $\psi \in V^*$ such that
$f(s) = \langle vs, \psi\rangle $. 

If $S=G$ is a group then $h(G) \subset \Is(V)$. 
\end{thm}
\begin{proof}
Apply Theorem \ref{main} to the flow $(S,S)$. Then for $f \in
\WAP(S)$ there exists a reflexive $V$ and a representation $h: S
\to \Theta(V)$, $\alpha : S \to B(V^*)$ such that $f(s)=\langle v,
\alpha(s)\rangle $ for a suitable $v \in V$. Denote by $e$ the identity
of $S$. Then
 $f= m_{v,\psi} $ where $\psi = \alpha(e)$.
\end{proof}

If we wish to get a {\it homomorphism}, just consider $h: S \to
\Theta(V)^{opp}= \Theta(V^*)$.

It is also easy now to establish the following result first established by
Shtern \cite{Sh} (see also \cite{Meop}).

\begin{f} \label{shtern}
The following are equivalent:
\begin{itemize}
\item
[(i)] A semitopological semigroup $S$ is wap (equivalently, $S$
is embedded into a compact semitopological monoid).
%(ii) There exists an $S$-wap flow $X$ with the strict action.
\item
[(ii)] There exists a reflexive space $E$ such that $S$ is
embedded (as a semitopological subsemigroup) into $\Theta(E)_w$.
%If $S=G$ is a topological group we can suppose that $G$ is embedded
%into $Is(V)_w=Is(V)_s$.
\end{itemize}
\end{f}
\begin{proof} (i) $\Longrightarrow$ (ii) We can suppose that $S$
is a monoid. Consider $X:=S^W$ the universal semitopological
compactification of $S$. Then the corresponding universal map
$u_W: S \to S^W$ is a topological embedding by (i) and hence, the
action $(S,S^W)$ is {\it left strict}. That is, there is no
strictly coarser topology on $S$ under which $S$ is a
semitopological semigroup and $S^W$ is still a semitopological
$S$-flow. By Theorem \ref{main} there exists a separating family
$(h_i, \alpha_i)$ of reflexive $V_i$-representations ($i\in I$) of
$(S,S^W)$. Then the $l_2$-sum of these representations defined on
the Banach space $V:=(\sum_iV_{i\in I})_{l_2}$ will induce a
weakly continuous antihomomorphism $h: S \to \Theta(V)$. Since the
original action is left strict, it is easy to show that $h$ must
be a topological embedding. Define $E:=V^*$. It is clear that the
antihomomorphism $h$ defines the desired homomorphism  $h: S \to
\Theta(V)^{opp}=\Theta(V^*)=\Theta(E)$.

\nt (ii) $\Longrightarrow$ (i) It is well known \cite{LG1} that
$\Theta(V)_w$ is a compact semitopological semigroup for every
reflexive $V$.
\end{proof}

By \cite{Meop}
(or, Corollary \ref{old} below), $\Is (V)_s=\Is (V)_w$ for every
reflexive $V$. Therefore we obtain the following result.

\begin{f} \cite{Meop}  Let $G$ be a topological group.
The following are equivalent:
\bit
\item
 [(i)] $G$ is wap.
\item
[(ii)] $G$ is a topological subgroup of the group $\Is(V)_s$
(endowed with the strong operator topology) of all linear
isometries for a suitable reflexive $V$. \eit
\end{f} 

\begin{remarks} \label{r-long} \ 

%(ii) We can suppose that $w(Z)=w(X)=w(Y)=d(Y)$.
%A $G$-flow $X$ is said to be
%{\it compactifiable} (or, $G$-Tychonoff \cite{Mesm})
%if there exists a $G$-compactification $\alpha: X \hookrightarrow Y$
%where $\alpha$ is a topological embedding into the compact jointly
%continuous $G$-flow $Y$.
\begin{itemize}
\item [(i)] Theorem \ref{wapeb} implies that
 every wap $S$-flow $X$ is compactifiable.
%\begin{corol} Every wap $(S,X)$ flow is compactifiable.
Moreover, if $S=G$ is a semitopological group and $\alpha: X \to
Y$ is a corresponding faithful wap $G$-compactification (which
exists by Theorem \ref{wapeb}) then the action $G\times Y \to Y$
is jointly continuous. This follows from Fact \ref{lawson}.
Therefore every noncompactifiable in a joint continuous way
$G$-space provides an example of a non-wap flow. Such examples can
be found even for jointly continuous
group actions of Polish topological groups $G$ on Polish
spaces $X$ (see \cite{Mesm, MeSc}).

%(iv) (Teleman)
%Every compactifiable $G$-flow can be embedded into the affine $G$-flow
%$(B_{V^*}, w^*)$ for some Banach space $V$.

\item  [(ii)] It is well known
(as noted for example in Arhangelskij \cite{Arh} or
Namioka-Wheeler \cite{NW}) that a compact space is Eberlein iff it
can be included into some right strict duality.
%This follows easily by Grothendieck's Lemma.
Theorem \ref{Davis} provides ``a flow version''.
\end{itemize}
\end{remarks}

%\sk

%Stability properties (EXTENT) of $[wap]^S$.

%\ssk
%subdirect products, factors for compact ...

%Hence, there exists a universal wap-compactification (Junghenn).
%Reformulate: wap iff the universal wap compactification is proper.

%\sk
%Corollary: $WAP_S(X)$ is a Banach algebra.
%(Hence (?), $\R$ or $\C$ is the same;
%$WAP(X,\R)$ separates iff  $WAP(X,\C)$.)

\section{Fragmentability and flows}

\vskip 0.1cm

\begin{defin} \label{d-fr}
Let $(X, \tau)$ be a topological space and $\mu$ be a uniformity
on the set $X$. Then $X$ is $(\tau,\mu)$-{\it fragmented} if for
each nonempty $A \subset X$ and each $\e \in \mu$ there exists a
$\tau$-open subset $O$ of $X$ such that $O\cap A \neq \emptyset$
and  $O\cap A$ is $\e$-small.
\end{defin}

This definition (for metrics) is explicitly defined by Jayne and
Rogers \cite{JR} and implicitly it appears even earlier in
Namioka-Phelps \cite{NP} (see also \cite{Ke}). There are several
generalizations: for covers (Bouziad \cite{Bo}), for functions
\cite{JOPV, Mefr}. Similar concepts are studied in many contexts: {\it
cliquish} (Thielman 1953),
 {\it huskable} (in French, {\it epluchable}) (Godefroy 1977).

The works \cite{Mefr,Meop} are devoted to a systematic study of the
fragmentability concept in the context of (semi)group actions and
topological dynamics.

Namioka's famous {\it joint continuity theorem} implies that every
weakly compact subset of a Banach space is norm-fragmented
\cite{Na2}. We need the following generalization for locally
convex spaces $(V, \mu)$ where $\mu$ denotes the usual additive uniform
structure on $V$.
\begin{lemma} \label{namiokamy} Every relatively
weakly compact $X \subset V$ in a l.c.s.
$V$ is $(weak, \mu)$-fragmented.
\end{lemma}
\begin{proof} See \cite[Proposition 3.5]{Mefr}.
\end{proof}

%\delete{
%\begin{defin} (Namioka pairs)
%$N(X,Y)$ means that $\forall$ sep. cont.
%$\Phi: X \times Y \rightarrow \R$ $\exists$ dense $G_{\delta}$-subset
%$D$ of $X$ s.t. $\Phi$ is j. cont. at $\forall$ point of $D\times Y$.
%\nl $X \in N$ if $N(X,Y)$ $\forall$ compact $Y$.
%\end{defin}
%
%\begin{lemma} \label{}
%Let X be a bounded subset of
%an l.c.s. $V$. If $N(X,K)$ $\forall$
%$w^*$-compact equic. $K \subset V^*$
%%(e.g. if $X \in N$)
%then $X$ is fragmented (w.r.t. $(\tau_w,\mu)$).
%\end{lemma}
%
%
%\begin{lemma} (Fragmentable subsets in l.c.s.)
%
%$X$ is fragmented in each of the following cases:
%
%\nt $\bullet$ $cl_w(X) \in N$.
%\nl $\bullet$ $cl_w(X)$ is Chech complete (e.g. w-compact
%$cl_w(fS)$
%
%($f\in WAP_S(Y)$, for instance)).
%\nl $\bullet$ $cl_w(X)$ is {\it $\alpha$-favorable} or
%{\it $\beta$-unfavorable}.
%\nl $\bullet$ $X$ is {\it dentable} in the sense of Rieffel
%
%(that is,
%$\forall \hssk U(0) \quad \exists a\in X: \hskip 0.3cm
%a \notin cl(co(X \setminus (a+U)).$
%\nl $\bullet$ $cl_w(X)$ is hereditarily Baire
%
%(see {\it huskable sets} in [Edgar-Wheeler 1984]).
%\end{lemma}
%}
%

We will use the following useful observation.

\begin{f} \label{Gdelta}
Let $(X,\tau)$ be a Baire space and $d$ a pseudometric on the
set $X$. If $X$ is $(\tau, d)$-fragmented, then $1_X: (X,\tau) \to
(X,d)$ is continuous at the points of a dense $G_{\delta}$ subset
$D$ of $X$.
\end{f}
\begin{proof}
Easily follows using the standard Baire arguments. See for example
the proof of Lemma 1.1 in \cite{Na2} or \cite[Lemma 3.2
(d)]{Mefr}.
\end{proof}
%This fact implies that every Eberlein compact has a dense $G_{\delta}$
%metrizable subset.

The following known characterization of Asplund spaces (which is a result of
many works) in terms of fragmentability is very important in our setting.

\begin{f} \label{asplundsp}
A Banach space $V$ is Asplund iff every bounded subset
$A \subset V^*$ of the dual $V^*$ is (weak$^*$,norm)-fragmented.
\end{f}

Standard examples of Asplund spaces include: reflexive spaces
and $c_0(\Gamma)$ spaces. Let $K$ be compact. Then $C(K) \in ASP$ iff
$K$ is {\it scattered} (that is, every
nonempty subspace of $K$ contains an isolated point).

Let $\mu$ be a uniformity on an $S$-flow $X$. We say:

\bit \item [a)] $z \in X$ is {\it a point of equicontinuity} (or,
a {\it Lyapunov stable}) (denote $z \in Equic_{\mu}(S,X)$ or,
simply, $z\in Equic$) if there exists a compatible uniformity
$\mu$ such that for all $\varepsilon >0$ there exists a
neighborhood $U(z)$ of $z$ such that $(sx,sz) \in \varepsilon$ for
every $(x,s)\in U\times S$.
\item [b)] $(S,X)$ is {\it (almost)
$\mu$-equicontinuous} if (resp.: $X=cl(Equic)$) $X=Equic$.

%\delete{
%(ii) $(S,X)$ is {\it orbitwise equicontinuous} if
%$z\in Equic(S,Sz)$ for every $z \in X$.
%
%(iii) {\it minimally equicontinuous} if every
%(not necessarily compact) minimal $S$-subspace is equicontinuous.
%
%Note that by our definition every minimally equicontinuous is orbitwise
%equicontinuous.
%}

\item
[c)] $(S,X)$ is  {\it uniformly $\mu$-equicontinuous} if for every
$\e \in \mu$ there exists $\delta \in \mu$ such that $(sx,sy) \in
\e $ for every $(x,y) \in \delta$ and every $s \in S$.

\item
[d)] A point $z\in X$ is the point of {\it local
$\mu$-equicontinuity} in the sense of Glasner and Weiss \cite
{GW2} if $z \in Equic_{\mu}(S, cl(Sz))$ (we do not require that
$X$ be compact). If this condition holds for every point in $X$, we
say that $(X, \mu)$ is locally equicontinuous and write $X
\in$ LE.

\item [e)] $(S,X)$ is  (almost, locally) equicontinuous
if $X$ is (resp.: almost, locally) $\mu$-equicontinuous with
respect to some compatible uniformity $\mu$ on $X$.

\item [f)] $(S,(X, \mu))$ is {\it not sensitive} (see for example
\cite{GW1} and the references there) if for every $\e \in \mu$
there exists a non-empty open subset $O$ of $X$ such that $(sx,sy)
\in \e $ for all $x,y \in O$ and $s \in S$.

\eit

\begin{thm} \label{Baire}  \ 
	\begin{itemize}
		\item [(i)] Every ${RN}$ (e.g., Eberlein) Baire flow $(S,X)$
		is almost equicontinuous. 
		\item [(ii)] 
		Every RN-approximable (e.g., wap) Polish $S$-flow
		$X$ is almost equicontinuous.
	\end{itemize}

\end{thm}

\begin{proof} (i)
There exists an Asplund representation $h: S \to \Theta(V)_w$,
$\alpha: X \to B(V^*)$. Then according to Fact
\ref{asplundsp}, $f(X)$ is (weak$^*$, norm)-fragmented.
The action of $\Theta(V)^{opp}$ on $(B(V^*), norm)$
is obviously uniformly equicontinuous.
Every point $z \in X$ of continuity of the map $1_X: (X,w^*) \to (X,norm)$
is a point of equicontinuity in the $S$-flow $(X,w^*)$.
Fact \ref{Gdelta} guarantees that such points are dense in $X$.
Therefore $(S,X)$ is almost equicontinuous.

(ii) Follows from (i) because every RN-approximable second
countable flow is RN (Corollary \ref{seccwap}).
\end{proof}

As a conclusion of the part (ii) and Corollary
\ref{c-wap-RN} we get the following known result.

%\delete{
%\begin{thm} \label{Polish}
%Every wap Polish  $S$-flow $X$ is almost equicontinuous.
%\end{thm}
%\begin{proof}
%Apply Theorem \ref{Baire} and Corollary \ref{countable} (iii).
%\end{proof}
%}

\begin{corol} \label{AAB}
(Akin-Auslander-Berg \cite{AAB}) Let $G$ be a topological group
and $X$ a metrizable compact $G$-flow. Assume that the $G$-flow $X$
is wap. Then $X$ is almost equicontinuous.
\end{corol}

\vskip 0.1cm

 The following definition is an important tool for our purposes.

\begin{defin} \label{defin-equifr}
Let $(X,\tau)$ be an $S$-flow and $\mu$ a uniformity on the set
$X$ such that $\tau \subset top(\mu)$. We say that the flow $(X,
\tau)$ is {\it $\mu$-equifragmented} if $X$ is
$(\tau,\mu)$-fragmented, the action of $S$ on $X$ is uniformly
$\mu$-equicontinuous and for some uniformity $\xi \subset \mu$ we
have $top(\xi)=\tau$.

%We say simply that {\it $X$ is equifragmented} if $X$ is
%$\mu$-equifragmented with respect to some $\mu$.
\end{defin}

%Let a monoid $S$ act on $X$ and $\psi$ be a uniformity on $X$.
%Denote by $\psi^{sup}$ the ``sup-uniformity'' on $X$ the base of
%which is generated by the system $\{\e^{sup}: \e \in \psi \}$,
%where
%$$\e^{sup} := \{ (x,y) \in X \times X : (sx,sy) \in \e \hskip 0.1cm \forall s
%\in S \}$$
%
% If $X$ is an $S$-flow with an admissible
%uniformity $\psi$ then $X$ is $\mu$-equifragmented with respect to
%some $\mu \supseteq \psi$ iff $X$ is $\psi^{sup}$-equifragmented.

We collect here some useful stability conditions for equifragmentability.

\begin{lemma}  \label{nearly} \ 
\bit
%\item [(i)] Let $X$ be $\mu$-equifragmented and $h: P\to
%X$ is continuous then the composition $f\circ h: P \to X$ is
%$\mu)$ is continuous.
\item
[(i)] The class of equifragmented flows is preserved under
subflows.
\item
[(ii)] Equifragmentability is preserved under products.
More precisely, if $X_i$
is $\mu_i$-eqifragmented then the product $ \prod X_i$ of
$S$-flows is $\prod \mu_i$-equifragmented.
\item [(iii)]
%\item [(iv)]
For every Asplund space $V$ the flow $(\Theta(V)^{opp}, (B(V^*), w^*))$
is $\mu_{\|\hskip 0.2cm  \|}$-equifragmented where
$\mu_{\|\hskip 0.2cm  \|}$ is the norm uniformity of $V^*$.
\item [(iv)] $(\Theta(V), (B(V), w))$ is $\mu_{\| \hskip 0.2cm
\|}$-equifragmented for every reflexive $V$.
\item [(v)] Every RN-approximable (e.g., wap) flow $(S,X)$ is
equifragmented.
%(with respect to a metric uniformity provided that
%$(S,X)$ is RN) .
\item [(vi)] If a compact flow $X$ is equifragmented then $X$ is not
sensitive. Therefore, every ${RN}$-approximable $S$-flow $X$ is
not sensitive.
 \eit
\end{lemma}

\begin{proof} The assertion (i) is trivial, (ii) and (vi) are
straightforward. For (iii) and (iv), we can use Fact
\ref{asplundsp} and Lemma \ref{namiokamy}, respectively. In order
to establish (v), combine (i), (ii) and (iii).
\end{proof}

Let a group $G$ act on a topological space $X$. We say that:

a) a point $z\in X$ is {\it transitive} (write: $z \in Trans$) if
$cl(Gz)=X$. If $Trans \neq \emptyset$, then, as usual,
$X$ is called transitive.

b) a point $z\in X$ is {\it quasitransitive}
(write: $z \in qTrans$) if $int(cl(Gz)) \neq \emptyset$.

c) $X$ is {\it quasiminimal} if $X=qTrans$.

d)  $X$ is {\it minimal} if $cl(Gz)=X$ for all $ z \in X$.

\vskip 0.1cm

%Let $(S, (X,\tau))$ be $\mu$-equifragmented.
 Let $(X,\tau)$ be a topological space and
 $\mu$ be a uniformity on $X$ such that
$\tau \subset top(\mu)$. We say that a subset $K \subset X$ is
$(\tau,\mu)$-{\it Kadec} if $\tau|_K=top(\mu)|_K$. Denote by
$Cont(\tau,\mu)$ the subset of all points of continuity of the
identity map $1_X: (X, \tau) \to (X, \mu)$. Clearly,
$Cont(\tau,\mu)$ is an example of a $(\tau,\mu)$-Kadec set.

\begin{thm} \label{genquasi}
Let a topologized group $G$ act on a topological space $(X,\tau)$
by homeomorphisms. If this action is $\mu$-equifragmented (with
respect to $\xi \subset \mu$ such that $top(\xi)=\tau$) then:

\bit \item [(i)] $qTrans  \subset Cont(\tau,\mu) \subset
Equic_{\mu}(G, X)$. In particular, every point of
quasitransitivity of $X$ is a point of $\xi$-equicontinuity.
\item
[(ii)] If a $G$-subflow $Y$ is quasiminimal (e.g., 1-orbit subset
$Y=Gz$) then $Y$ is a $(\tau,\mu)$-Kadec set. Hence $\mu|_Y$ is an
compatible uniformity on $Y$ and $Y$ is a uniformly
$\mu|_Y$-equicontinuous $G$-flow. \eit
\end{thm}

\begin{proof} (i) Let $z \in qTrans$.
We have to show that for every $\varepsilon \in \mu$ there exists
a $\tau$-neighborhood $O(z)$ of $z$ such that $O$ is
$\varepsilon$-small. Choose $\delta \in \mu$ such that
$(gy_1,gy_2) \in \varepsilon$ for every $(y_1,y_2) \in \delta$ and
$g \in G$. Since $z \in qTrans$, the set $A:=int(cl(Gz))$ is
non-void. Since $ X$ is  $(\tau,\mu)$-fragmentable, we can pick a
{\it non-void} $\tau$-open subset $W$ of $X$ such that $W\subset
A $ and $W$ is $\delta$-small in $X$. Clearly, $W \cap Gz \neq
\emptyset$. One can choose $g_0 \in G$ such that $g_0 z \in W \cap
Gz$. Denote by $O$ the open subset $g_0^{-1} W$ of $(X, \tau)$.
Then $O$ is a  $\tau$-neighborhood of $z$ and is
$\varepsilon$-small. This proves the inclusion
$qTrans  \subset Cont(\tau,\mu)$. The second inclusion
$Cont(\tau,\mu) \subset Equic_{\mu}(G, X)$ is trivial because $X$ is
uniformly $\mu$-equicontinuous.

(ii) By the quasiminimality of $Y$, $qTrans(Y)=Y$. Therefore, the
assertion (i) implies that $\tau|_Y=top(\mu)|_Y$.
\end{proof}

\begin{thm} \label{asplundr}
Let $G$ be a semitopological group and $X$ be an RN-approximable
semitopological $G$-flow.
Then $X \in LE$ and every $G$-quasiminimal subspace (for instance,
every orbit) of $X$ is equicontinuous.
\end{thm}
\begin{proof} By  Lemma \ref{nearly} (v), $X$ is $\mu$-equifragmented.
For every fixed $z \in X$ consider the $S$-subflow $Y:=cl(Gz)$.
 Clearly, $z$ is a point of quasitransitivity of $Y$.
Then we can apply Theorem \ref{genquasi} to $(G,Y)$ and conclude that
$z$ is a point of local equicontinuity of $Y$ (and hence of $X$).
\end{proof}

\begin{corol} \label{gVTA}
(``Generalized Veech-Troallic-Auslander Theorem'')

\nt Every wap (not necessarily compact or metrizable) $G$-flow $X$ is
LE and every $G$-quasiminimal subflow of $X$ is equicontinuous.
\end{corol}
\begin{proof} By Corollary \ref{c-wap-RN}
every wap flow is RN-approximable.
\end{proof}
%$\mu_{wap}$-$LE$.
%Such $X$ is minimally equicontinuous (and, hence,
%orbitwise equicontinuous).

%\delete{
%\begin{proof}
%Denote by $\mu_{wap}$ the precompact uniformity generated by
%$WAP_S(X)$ on a wap $S$-space $X$. Then by
%Theorems \ref{genquasi} and \ref{nearly} $X$ is LE with respect to
%the {\it admissible} uniformity $\mu_{wap}$ .
%\end{proof}
%}

\begin{remark} \label{r-VTA}
Troallic \cite{Tr} and also  Auslander \cite{Au} proved
that every minimal compact wap $G$-flow $X$ is equicontinuous.
Previously such a result was established for compact Eberlein (in our
terminology) $G$-flows by Veech \cite{Ve}.
%About the inclusion: ``WAP
%$\subset$ LE''; it was known only that every compact
%metrizable wap $G$-flow is LE (see \cite{GW2}).
\end{remark}

Combining Theorem \ref{wapeb} and Corollary \ref{gVTA}
for general $G$-flows, we can draw the following diagram
%\begin{diagram}[height=1.7em,width=2.3em,tight]
%  Eb   & \rTo & RN  & \rTo & \{RN-approximable\}
%& \rTo & LE  \\
%            & \rdTo &  & &  &\ruTo{}   \\
%                                       &    & {WAP}  \\
%\end{diagram}
$$\text{Eberlein} \subset \text{WAP} = \text{REFL}_{app}
\subset \text{RN}_{app} \subset \text{LE}$$
%For second countable $G$-flows we have
%$$Eberlein = wap = REFL-approximable \to RN = RN-approximable \to LE$$
%Clearly for every fixed $S$ holds $WAP \subset RN \subset LE$ and
%$RN \subset AE$.
%The inclusions are proper.

Consider the case of $S=\Z$ and metrizable compact cascades. The
fact that WAP $\neq$ LE is discussed in \cite{GW2}. The authors
constructed (see main example in \cite[page 350]{GW2}) a
transitive cascade $(\Z, X)$ such that $X$ is in LE but not wap
and every point of transitivity is {\it recurrent}. If we do not
require the last assumption then there exists an elementary example
distinguishing wap and RN (and , hence also wap and LE). Namely,
the two-point compactification $X$ of $\Z$  with the
natural action of $\Z$ on $X$ is transitive and
contains two fixed points. Fact \ref{unique} implies that
such $(Z,X)$ can not be wap. On the other hand, $X$ is clearly scattered.
Therefore $(Z,X)$ is RN by Proposition \ref{scattered} below.

%\delete{
%In the next example we give a simpler cascade which distinguishes
%$WAP$ and $RN$. In particular, it distinguishes  $LE$ and $WAP$.
%Example. Choose a point $0<t_0 <1$
%and define an$K\subset V$ is a w-compact S-invariant via
%$h: S \to \Theta (V)$.
% orbit closure
%$X=\{0\}\cup \{1\} \cup \{a_n=\sigma^n (t_0) \ n \in \Z \}$. Then we have
%a subflow of an example abow. Recall that $\sigma(x)=x^2$.
%Each non-fixed point $a_n$ is isolated in $X$. Hence,
%$X$ is locally equicontinuous, $X \in LE$. Moreover, $X$ as a $\Z$-flow
%is $RN$. Indeed, $X$ is scattered hence as it is well known
%(see Zemadeni) $C(X)$ is Asplund. By Stegall's theorem, $C(X)^*$ has RNP.
%Now recall that $X$ is a $Z$-subflow of $B^*_{w^*}$. Therefore,
%$(\Z, X) \in RN$. On the other hand,
%$(\Z, X) \notin WAP$ because Lemma below ($X$ is transitive and
%contains two minimal subsets $\{0\}$ and $\{1\}$).
%}

\begin{f}
\label{unique}
Let $X$ be a wap transitive compact $G$-flow.
Then $X$ contains a unique minimal compact subflow.
\end{f}
\begin{proof}
Let $E=E(G,X)$ be the Ellis (semitopological) semigroup.
By \cite[Proposition II.5]{EN} this semigroup contains a
unique minimal ideal $K$ which is closed in $E$.
It follows by transitivity that $Et_0=X$ for some $t_0 \in X$.
Then the unique minimal compact subset of $X$ is $Kt_0$.
\end{proof}
%
%For example, the cascade $(\Z,[0,1])$ generated by the map
%$f(x)=x^2$ on the closed interval  is not wap.

\delete{
\begin{example} \label{skew}
Define a distal cascade $(\Z, X)$ on the 2-dimensional torus $X=\T^2$ by
$$\Z \times \T^2 \to \T^2, \hskip 2cm n \circ (a,b)=(a+n\cdot b,b).$$
Then $X \in$ LE but $X$ has no points of equicontinuity.
In particular, such $(\Z, X)$ is not RN by Theorem \ref{Baire}.
\end{example}
}

If the group action $(G,X)$ is RN then there exists a compatible
uniformity $\mu$ on $X$ (the precompact uniformity of the
corresponding weak star $G$-embedding of $X$ into $B(V^*)$ with
Asplund $V$) such that $X$ is not sensitive (see Lemma
\ref{nearly} (vi)). Another observation comes from Theorem
\ref{asplundr}. It implies that every RN-approximable
1-orbit group action is
equicontinuos. This provides an easy way producing examples of
$G$-flows which fail to be RN. Roughly speaking, RN $G$-flow
cannot be ``too chaotic'' or ``too massive''.

%For example 1-orbit ``massive actions'' (like
%$\approx$ 2-transitive are not RN).
%More precisely, suppose that  $(G,X)$ is a 1-orbit flow and
%for a certain pair $a\neq b, z \in X$ and a net $g_i \in G$
%holds (the topological proximality)
%$g_ia \rightarrow z \leftarrow g_ib$ then $(G,X)$ is not RN.

%NIL-EXAMPLE

\vskip 0.3cm

Let $G$ be a topological group. Consider the natural action
(call it a ``$\Delta$-action'')
$$\pi_{\Delta}: (G\times G) \times G \to G, \quad (s,t)x=sxt^{-1}.$$
%It actually coincides with the coset $G$-space action $(G\times G,
%X / H)$, where $H=\Delta:= \{(g,g) : g \in G\}$.

\begin{lemma} \label{delta}
Let $G$ be a topological group such that
$(G\times G, G, \pi_{\Delta})$ is an
RN-approximable (e.g., wap) flow. Then
$G$ satisfies SIN
%%may29 PROPERTY
({\it small invariant neighborhoods}).
%(by Theorem \ref{asplundr}).
\end{lemma}
\begin{proof}
The given (1-orbit)
$\Delta$-action is $\mu$-equicontinuous, by Theorem \ref{asplundr},
with respect to some compatible uniformity $\mu$ on $G$. Let
$U(e)$ be an arbitrary neighborhood of the identity in $G$. Choose
$\varepsilon \in \mu$ such that the neighborhood $\e (e)=\{x \in
G: (e,x) \in \e \}$ is contained in $U(e)$. By the
$\mu$-equicontinuity of the $\Delta$-action at the point $e$ one
can choose a neighborhood $O(e)$ such that $ sOt^{-1}$ is
$\varepsilon$-small for all $(s,t) \in G \times G$. Then $gOg^{-1}
\subset \varepsilon(e) \subset U(e)$ for every $g \in G$.
This is equivalent to the condition $G\in SIN$.
%({\it small invariant
%neighborhoods}).
%(for the definition and properties see \cite{}).
\end{proof}

%%This lemma provides many not RN 1-orbit group actions
%%with locally compact (or, even Lie) acting groups.

%\delete{
%Recall that by [Uspenskij 1993]
%$St_z \hookrightarrow H(S_1)$ is not an {\it epic}.
%
%\begin{thm}
%Let $H$ be a proper subgroup of $G$ and the coset $(G,G/H)$ space
%is RN. Then  the inclusion $H \hookrightarrow G$ is not {\it epic}.
%\end{thm}
%\begin{proof} By Theorem \ref{1}
%$\quad (G,G/H)$ is equicontinuous. Then by \cite{MeF},
%$G/H$ is $G$-linearizable. Therefore,
%{\it Free topological $G$-group} on $G/H$ is not trivial.
%By Pestov's result (see \cite{}) $H \hookrightarrow G$ is not epic.
%\end{proof}
%}

\vskip 0.1cm Now we can strengthen a  result of Hansel and
Troallic. Let $G$ be a topological group. Following \cite{HT} we
say that a function $f\in C(G)$ is {\it strictly wap} (notation:
$f \in s\WAP(G)$) if $GfG$ is relatively w-compact in $C(G)$.
Denote by
 $[swap]$ the class of groups such that $s\WAP(G)$
separates points and and closed subsets of $G$. Denote by $[WS]$
the class of groups for which every wap function is strictly wap.
Clearly, $[wap] \cap [WS] \subset [swap]$.

%Define also $sWAP(G)$ $=\{$ strictly wap functions on $G \}$

\begin{prop}
If  $G \in [swap]$ then $G \in SIN$.
\end{prop}
\begin{proof} First observe that
$ G \in [swap]$ iff the $\Delta$-action
$(G\times G, G)$ is wap and hence
RN-approximable.
Now Lemma \ref{delta} finishes the proof.
\end{proof}

\begin{corol} (Hansel-Troallic \cite{HT})
Let $G\in [wap]\cap [WS]$. Then $G \in SIN$.
\end{corol}

%\delete{
%\sk
%
%\section{Questions}
%
%\sk
%
%\begin{question}
%
%\end{question}
%}

\section{Asplund functions and representations} \label{s-asp-f}

\vskip 0.1cm

Recall that a Radon-Nikod\'ym compact space \cite{Na1} is a
compact subset in $(V^*, w^*)$ for an Asplund space $V$. We
introduce a generalization for flows. Our approach synthesizes
some ideas from \cite{Ve, St, Na2, Fa}.
%For instance the proof of Theorem
%\ref{thm1} is based on the results and ideas of \cite{Fa}.

The following definition goes back to Stegall \cite{St} and
Namioka \cite {Na2}.

\begin{defin} \label{d-asplset}
Let $M$ be a nonempty bounded subset of a Banach space $V$.
Say that $M$ is an {\it Asplund set} in $V$ if for every countable $C
\subset M$ the pseudometric space $(V^*,\rho_C)$ is separable,
where $$\rho_C(\xi, \eta)= sup \{|\langle c, \xi\rangle  - \langle c, \eta\rangle |: \hskip 0.2cm c
\in C \}.$$
\end{defin}
%We use the notation $M \in A(K)$.

Generalizing slightly this
definition, let's say that $M$ is an {\it Asplund set for} $K \subset V^*$
if the pseudometric subspace $(K, \rho_C)$ is separable for every
countable $C \subset M$.

We need the following lemma of Namioka in the form presented by
Fabian.

\begin{lemma} \label{Na-Fa}
\cite[Lemma 1.5.3]{Fa} Let $X$ be a compact space (canonically
embedded into $C(X)^*$) and let $M \subset C(X)$ be a bounded
subset. Assume that $(X, \rho_M)$ is separable. Then the
pseudometric space $(C(X)^*, \rho_M)$ is also separable.
\end{lemma}

\begin{corol} \label{c:Na-Fa}
$M \subset C(X)$ is an Asplund set for compact $X$ iff $M$ is an
Asplund set for $C(X)^*$.
\end{corol}

\begin{remark} \label{r-asp-pr}
The family of Asplund sets in $V$ has  nice properties being
stable under taking subsets, finite unions,
closures, linear continuous images,
finite linear combinations, etc. Note also that if $M_1$ and $M_2$
are Asplund sets in $C(X)$ for a compact $X$, then the subset $M_1
\cdot M_2$ is also Asplund. For these and some other results we
refer to \cite{Fa}.
\end{remark}

We say that a {\it bounded}
duality $Y \times X \to \R$ is an {\it Asplund
duality} if $q_Y(Y)$ is an Asplund subset of $C(X)$. Conversely,
the subset $M \subset C(X)$ is an Asplund set iff the
corresponding duality $M \times X \to \R$ is an Asplund duality
where $M$ is endowed with the pointwise topology inherited from
$C_p(X)$.

The following Lemma is a reformulation of a result of
Namioka \cite[Theorem 3.4]{Na2}.

\begin{lemma} \label{l-aspl-du}
Let $V$ be a Banach space and $K$ be a compact subspace in the
dual ball $B^*_{w^*}$. Suppose that $K$ is (weak$^*$,
norm)-fragmented in $V^*$. Then $B=B(V)$ is an Asplund set for $K$
and $B \times K \to \R$ is an Asplund duality for every topology
$\tau$ on $B$ such that $\psi: B \to \R$ is $\tau$-continuous for
every $\psi \in K$.
\end{lemma}

By Fact \ref{asplundsp} and Lemma \ref{l-aspl-du}, the strict duality
$rB_w \times B^*_{w^*} \to \R$ is an Asplund duality
(call it: {\it canonical Asplund duality}) for every Asplund space $V$ (and $r>0$).

\begin{defin} \label{def-w-admis}
Let $X$ be a compact $S$-flow with a separately continuous left
action. We say that $X$ is {\it $w$-admissible } if
$C(X)=\WRUC_S(X)$.
\end{defin}

This happens, for example, if either: a) $(S,X)$ is jointly continuous (then by
Fact \ref{f-c=r} we have even  $C(X)=\RUC_S(X)$); b) $(S,X)$ is wap
(by Fact \ref{f-functions}); or c) $S$ is a $k$-space (use
Fact \ref{Groth}).

\begin{thm} \label{thm-aspldual}
Let $X$ be a compact $w$-admissible $S$-flow. Every Asplund strict
%mich
$S$-duality $Y \times X \to [-r,r]$ is an $S$-restriction of a
canonical Asplund $S$-duality
%$B \times B^* \rightarrow \R$
with respect to a weakly continuous antihomomorphism $h: S \to
\Theta(V)$. More precisely, there exist: a suitable Asplund space
$V$ and equivariant maps: $\gamma_1: X \to B^*_{w^*}$ and
$\gamma_2: Y \to B$ such that the following diagram is
commutative

\begin{equation*} \label{diag4}
\xymatrix{
Y \ar@<-2ex>[d]_{\gamma_2} \times X
\ar@<2ex>[d]^{\gamma_1} \ar[r]  & [-r,r] \ar[d]^{id} \\
rB \times B^* \ar[r]  &  [-r,r] }
\end{equation*}
  
%nd
%\begin{diagram}[height=2.1em,width=3.0em,tight]
% &Y \times X &\rTo^{} & \R  \\
% &\dTo^{\gamma_2} \quad \quad \dTo_{\gamma_1}  &    &\dTo_{id}  \\
% &rB \times B^* & \rTo^{} & \R \\
%\end{diagram}
where we require that $\gamma_1$ is a topological embedding and
$\gamma_2$ is an injective map.

 If the action $S \times X \to X$ is
jointly continuous, then we can suppose that $h: S \to \Theta(V)$
is strongly continuous.
 \end{thm}
\begin{proof}
% We can suppose that $S$ is a monoid and $ex=x$.
It is enough to prove for $r=1$. 
 Consider the natural continuous injective map $q_Y: Y \to C_p(X)$ and
 denote by $K$ the subset $q_Y(Y) \subset C(X)$. Then $K$ is an Asplund
subset in $C(X)$. Moreover, $K$ is a subset in the unit ball of $C(X)$.
 Then $K$ is an Asplund subset also in the Banach subspace $E=cl(sp(K))$ of $C(X)$.
Following the method of \cite{St} and, especially,
\cite[Section 1.4]{Fa}, one can modify the proof
of Theorem \ref{Davis} using the factorization procedure for
Asplund $S$-sets (instead of weakly compact sets). We define a
sequence $\| \cdot \|_n$ of norms on $E$ each of them equivalent
to the original norm. Namely, for every natural $n$ consider Minkowski's
functional of the set 
$$P_n:= a^{\frac{n}{2}} co(-K \cup K) + a^{-\frac{n}{2}} B(E),$$
where, as in the proof of Theorem \ref{Davis}, 
$a$ is the unique solution of the equation $\sum_{n=1}^{\infty} \frac{a^n}{(a^n+1)^2} =1$

It is important that the subset $\cap_{n \in \N} P_n$ is Asplund.
Moreover, by \cite[Theorem 1.4.4]{Fa} we get a linear injective
continuous mapping $j: V \to E$
%mich
%(with $\|j\| \leq 1$)
where $V$ is an Asplund space. Since $K$ is an $S$-invariant subset of $C(X)$,
the same is true for $E$, $B(E)$ and $P_n$. Therefore, every norm
$||\cdot ||_n$ is $S$-nonexpansive. Then it follows by the
construction that the corresponding norm $N$ on the Banach space
$(V, N)$ is also $S$-nonexpansive.

Define $\gamma_2: Y \to B(V)$ as a natural $S$-inclusion (of
sets). Since $X$ is $w$-admissible we have $\WRUC_S(X)=C(X)$. This
guarantees that every orbit map ${\tilde z}: S \to C(X)$ is weakly
continuous. Hence the action of $S$ on $(E, weak)$ is separately
continuous.

On the other hand, by \cite[Theorem 1.4.4]{Fa}, the adjoint map
$j^*: E^* \to V^*$ has the norm dense range. It follows that for
every bounded subset $A$ of $V$, the weak topology of $V$ and the
weak topology of $E$, considering of $A$ as a subset of $E$ and
$C(X)$, are the same. In particular, this implies that every orbit
map ${\tilde v}: S \to (V,w)$ is weakly continuous. Thus, the
antihomomorphism $h: S \to \Theta(V)$ is weakly continuous.

We get also that the dual (left) action of $S$ on $V^*$ is weak*
separately continuous. The natural $S$-inclusion $j: V \to C(X)$
is a linear continuous $S$-map.
% with a dense range $j(V)$ because
%$Q \subset j(V)$ and $cl(sp(Q))=C(X)$. Therefore, the adjoint
%$j^*: C(X)^* \to V^*$ is an injective $S$-map. Since $j^*$ is
%weak$^*$ weak$^*$ continuous we obtain that $B(C(X)^*)$ is
%canonically $S$-homeomorphic to a weak$^*$ compact $S$-subset of
%$V^*$ where $V$ is Asplund.
The adjoint $j^*: C(X)^* \to V^*$ is a weak$^*$-weak$^*$
continuous $S$-operator. Denote by $\gamma_1$ the restriction of
this map on $X \subset C(X)^*$. Clearly, $\langle y,x\rangle =\langle \gamma_2(y),
\gamma_1(x)\rangle $. Then $\gamma_1$ is injective (and hence a
topological embedding) because the original duality is (right)
strict.

If $(S,X)$ is a jointly continuous flow then, like Theorem
\ref{Davis}, we can prove that $h$ is strongly
continuous, too. 
%use that $C(X)=RUC_S(X)$.
\end{proof}

 It is well known (see \cite{Na2, Re}) that,
similarly to the ``Eberlein case'',
a compact space $X$ is RN iff the unit ball $B^*\subset C(X)^* $
(and hence P(X)) is RN. The following result provides, in particular,
a generalization for flows.

\begin{thm} \label{thm1}
Let $(S,X)$ be a compact $w$-admissible  (e.g., jointly
continuous) flow. The following conditions are equivalent:
\begin{enumerate}
\item [(i)]
The flow $(S,X)$ is RN.
\item [(ii)] There exists a representation $(h, \alpha)$ of $(S,X)$ into a
Banach space $V$ such that the antihomomorphism
$h: S \to \Theta(V)$ is weakly
continuous, $\alpha : X \to V^*$ is a bounded
weak$^*$ embedding, and $\alpha(X)$ is (weak$^*$, norm)-fragmented.
\item [(iii)] There exists a representation $(h, \alpha)$ of $(S,X)$ into
a Banach space $V$ such that $h:~S \to \Theta(V)$ is an
antihomomorphism (no continuity assumptions on $h$), $\alpha : X \to V^*$
is a bounded weak$^*$ embedding, and $\alpha(X)$ is (weak$^*$,
norm)-fragmented.
\item [(iv)] There exists a (right) strict
Asplund $S$-duality $Y \times X \to \R$.
\item [(v)]
There exists a bounded $S$-invariant subset $M \subset C(X)$ such
that $M$ separates points of $X$ and $M$ is an Asplund set for $X$
(equivalently, for $(C(X)^*)$).
\item [(vi)]
There exists an $S$-invariant Asplund set $Q$ in $C(X)$ such that
cl(sp(Q))=C(X).
\item [(vii)]
$(S, B(C(X)^*)$ is RN.
\item [(viii)]
$(S,P(X))$ is RN.
\end{enumerate}

If $(S,X)$ is jointly continuous, then in the assertions (i), (ii),
(vii), (viii) we can suppose in addition that the corresponding
$h: S \to \Theta(V)$ is strongly continuous.
\end{thm}
\begin{proof}

(i) $\Longrightarrow$ (ii) By definition of RN, there exists a faithful
Asplund $V$-representation. By Fact \ref{asplundsp},
$\a(X)$ is $(w^*, norm)$-fragmented in $V^*$.

(iii) $\Longrightarrow$ (iv) The ball $B(V)$ is $S$-invariant and
separates points of $\alpha(X)$. It follows from Lemma
\ref{l-aspl-du} that the right strict $S$-duality
$$B(V) \times X \to [-1,1], \hskip 0.1cm (v,x) \mapsto \langle v, \a (x)\rangle $$
is an Asplund duality. In order to get a strict duality, pass to
the associated reduced duality $\langle ,\rangle_q: B_q \times X \to [-1,1]$
(using the ``dual version'' of Lemma
\ref{net}.4). Clearly, $\langle ,\rangle_q$ is also an Asplund duality.

% Consider the introversion type operator $T: V \to C(X)$,
%$T(v)(x)=<v, \alpha(x)>$. Then $M:=T(B(V))$ is a bounded
%$S$-invariant Asplund set for $X$. By Lemma \ref{Na-Fa}, $M$ is an
%Asplund set (for $C(X)^*$) in the Banach space $C(X)$. The family
%$M=\{T(v): v \in B(V) \}$ separates points of $X$. Hence,

(iv) $\Longrightarrow$ (v) Take $M:=q_Y(Y) \subset C(X)$
(and use Corollary \ref{c:Na-Fa}).

(v) $\Longrightarrow$ (vi) Suppose that a set $ M$ satisfies
assumptions of (v). As in the proof of Theorem \ref{AL}, produce
inductively the sequence of subsets $M_n=M_1 \cdot M_1 \cdots
M_1$, where $M_1=M \cup \{1\}$. Then it is well known that each
$M_n$ is again an Asplund set (Remark \ref{r-asp-pr}). Moreover,
even the set $Q = \cup_{n\in \N} 2^{-n} M_n$ is Asplund. On the
other hand, by the Stone-Weierstrass theorem $cl(sp(Q))=C(X) $. By the
construction $Q$ is $S$-invariant.
%Thus, $C(X)$ is Asplund $S$-generated.

(vi) $\Longrightarrow$ (vii) Since $Q$ is an Asplund set in $C(X)$
we obtain that $Q \times B^* \to \R$ is a (left strict) Asplund
$S$-duality. This duality actually is right strict (and hence
strict) because $cl(sp(Q))=C(X)$. Now we can conclude that $B^*$
is RN $S$-flow as it follows directly from Theorem
\ref{thm-aspldual}. The same result guarantees that
$h$ is strongly continuous provided that
$(S,X)$ is a jointly continuous flow.

Other implications are trivial.
%If $(S,X)$ is a jointly continuous flow then like Theorem
%\ref{davis} we can use that $C(X)=RUC_S(X)$.
\end{proof}

\begin{defin} \label{d-ASPf}
 Let $X$ be a semitopological compact
$S$-flow. Let us say that a function $f \in C(X)$
%(just $f \in C(X)$, if the action is
%$w$-admissible)
is {\it $S$-Asplund} if the orbit $fS$ is an Asplund set in
$C(X)$. Equivalently, if $(fS) \times X \to \R$ is an Asplund
duality. More explicitly, taking into account Corollary
\ref{c:Na-Fa}, we see that $f$ is Asplund iff for every countable
subset (equivalently, separable) $C \subset S$ the pseudometric
space $(X,\rho_C)$ is separable, where
$$\rho_C(x, y)= sup \{|f(sx) - f(sy)|:
\hskip 0.1cm s \in C \}.$$ If $S$ is separable then it is
equivalent to check the separability of the single semimetric
space $(X,\rho_S)$.
\end{defin}

Denote by $\Asp_S(X)$ the set of all $S$-Asplund functions on a compact $X$.
%By the definition $Asp_S(X) \subset WRUC_S(X)$.
The product $F=f_1 f_2$ of two Asplund functions $f_1$ and
$f_2$ on $X$ is again Asplund. This follows from
the inclusion $FG \subset (f_1G)\cdot (f_2G)$ taking into
account Remark \ref{r-asp-pr}. It is easy to show
that in fact $\Asp_S(X)$ is a Banach $S$-subalgebra of $C(X)$ for
every
%!!!!!!!!! $w$-admissible (e.g., jointly continuous)
compact $S$-flow $X$.

%\delete{
%Let $i_R: S \to S^R$ be
%the universal jointly continuous $S$-compactification of the left action
%$(S,S)$.
%Denote by  $Asp(S)$ the set of all functions $f: S \to \R$ such that
%for some Asplund function $F: S^R \to \R$ on the jointly continuous compact
%$S$-flow $S^R$ holds $f(s)=F(i_R(s))$ for every $s \in S$.
%Thus, by the definition $Asp(S) \subset RUC(S)$.
%}

\begin{lemma} \label{lemma-functions}
$\WAP_S(X) \subset \Asp_S(X)$
for every semitopological compact $S$-flow $X$.
\end{lemma}
\begin{proof}
Every weakly compact subset of a Banach space is
an Asplund set (see \cite{Fa}). In particular, $fS \subset C(X)$
is an Asplund set if $fS$ is relatively weakly compact.
\end{proof}

\begin{thm} \label{RNrepres}
({\bf RN Representation Theorem}) Let $X$ be a compact
$w$-admissible $S$-flow. The following are equivalent:
\begin{enumerate}
\item [(i)] $f \in \Asp_S(X)$.
\item [(ii)]
$(fS) \times X \to \R$ is an Asplund $S$-duality.
\item [(iii)]
There exist: a representation $(h,\alpha)$ of  $(S,X)$ into an
Asplund $V$ with a (weakly continuous) antihomomorphism $h$ and a
vector $v \in V$ such that $f(x) =\langle v, \alpha(x)\rangle $.
\item [(iv)]
There exist: a representation $(h,\alpha)$ of $(S,X)$ into a Banach
space $V$ with a weak$^*$
%mm22Jan07    weak$^*$-continuous
$\a: X \to B(V^*)$ (no continuity assumptions on $h$) such that
$\a(X)$ is $(w^*, norm)$-fragmented and there exists a vector $v
\in V$ satisfying $f(x)=\langle v, \a(x)\rangle $ for every $x \in X$.

\end{enumerate}
If $(S,X)$ is jointly continuous, ``weakly continuous'' can be replaced by
``strongly continuous''.
\end{thm}
\begin{proof}
(i) $\Longleftrightarrow$ (ii) Is trivial as it was mentioned
earlier.

(ii) $\Longrightarrow$ (iii) Using Lemma \ref{net}.4
pass to the associated strict $S$-duality
$(fS) \times X_q \to \R$ (which again is
Asplund) and apply Theorem \ref{thm-aspldual}.

(iii) $\Longrightarrow$ (iv) is trivial by Fact \ref{asplundsp}.

(iv) $\Longrightarrow$ (i) Observe that the orbit $vS$ is an
Asplund set for $\alpha(X)$ by Lemma \ref{l-aspl-du}. Now observe
that the orbit $T(vS)=fS$ of $f$ is an Asplund set for $X$ (and
hence for $C(X)^*$ by Corollary \ref{c:Na-Fa}).
%Observe that for every $x \in X$ we have
%$(fs)(x)=f(sx)=<v, \alpha(sx)>$. That is,
%$fS$ behaves like $Sv$ on $X$. On the other hand , $Sv$ is an Asplund set
%in $V$ for $\alpha(X)$ by \cite[Theorem 3.4]{Na2}.
%%
%%
%(i) $\Longrightarrow$ (ii) Consider the continuous map $i_f: X \to
%C_p(fS)$, $i_f(x)(fs)=f(sx)$ and denote by $X_f$ the compact
%space $i_f(X)$. On the other hand consider also the map $j_f: fS
%\to C(X_f)$, $j_f(fs)(x)=f(sx)$. Then it is easy to see that the
%subset $M:=j_f(fS) \subset C(X_f) \subset C(X)$ is an Asplund set
%for $X_f$.
% By Theorem \ref{thm1} there exists an
% Asplund representation $(h, \alpha_f)$ in $V$ with a weakly
%continuous $h: S \to \Theta(V)$ and a weakly$^*$ continuous {\it
%embedding} $\alpha_f : X_f \to B(V^*) $. By the construction $M$
%is a subset of $V$ (see the proof of Theorem \ref{thm1}, part (v)
%$\Rightarrow$ (vi), where the constructed Asplund space is denoted
%by $Z$) . Now for the desired $v$,choose $v=j_f(f) \in V $ and
%define $\alpha: X \to B(V^*)$ as the composition of two maps:
%$i_f: X \to X_f$ and $\alpha_f$.
\end{proof}

\begin{corol} \label{compasp}
Let $X$ be a compact $w$-admissible $S$-flow.
The following are equivalent:

(i) $X$ is RN-approximable $S$-flow.

(ii) $C(X)=\Asp_S(X)$.
\end{corol}
\begin{proof}
(i) $\Longrightarrow$ (ii) Let $X$ be RN-approximable. Then by
Theorem \ref{RNrepres}, $\Asp_S(X)$ separates points of $X$. On the
other hand, $\Asp_S(X)$ is a Banach subalgebra of $C(X)$ containing
the constants. Therefore by the Stone-Weierstrass Theorem we have
the coincidence $\Asp_S(X)=C(X)$.

(ii) $\Longrightarrow$ (i) The algebra $C(X)=\Asp_S(X)$ separates
points and closed subsets of $X$. Hence, by Theorem \ref{RNrepres}
there are sufficiently many Asplund representations of $(S,X)$.
\end{proof}

\begin{prop} \label{p-factorRN} \ 
\ben \item Let $X$ be a compact $S$-flow and $q: X \to Y$ be an
$S$-quotient. A continuous bounded function $f: Y \to \R$ is
$S$-Asplund iff the composition $F=f\circ q : X \to \R$ is
$S$-Asplund.
 \item $F \in \Asp(X)$ iff it comes from an RN
$S$-factor. That is, there exist an RN $S$-flow $Y$, an $S$-factor
$q: X \to Y$ and a continuous function $f \in C(Y)$ such that $F=
f \circ q$. \item A factor $Y$ of RN-approximable compact
$w$-admissible $S$-flow $X$ is again $RN$-approximable and
$w$-admissible. If $Y$ is metrizable then $(S,Y)$ is an RN flow.

\een
\end{prop}
\begin{proof} (1) For every pair $x_1,x_2$ in $X$
and every countable $C \subset G$ we have
$$\rho_{FC} (x_1,x_2)=sup \{ |F(sx_1)-F(sx_2)|: \hskip 0.1cm s \in C \}=$$
$$=sup\{|f(sq(x_1))-f(sq(x_2))|: \hskip 0.1cm s \in C\}=\rho_{fC}
(q(x_1),q(x_2)).$$

 Thus, $q: (X, \rho_{FC}) \to (Y,\rho_{fC})$ is a surjective
 ``pseudometric-preserving'' map.
In particular,
 $(X,\rho_{FC})$ is separable iff $(Y,\rho_{fC})$ is separable.

(2) Combine the first assertion, Corollary \ref{compasp} and
RN representation Theorem \ref{RNrepres}.

(3) 
First observe that $f\in \WRUC_S(Y)$ iff $F \in \WRUC_S(X)$
because $q: X \to Y$ induces the $S$-inclusion
$q^*: C(Y) \hookrightarrow C(X)$ of Banach $S$-algebras.
This implies that $Y$ is also $w$-admissible.
We have to show that $Y$ is RN-approximable.
By Theorem \ref{RNrepres} it is
equivalent to check that $C(Y)=\Asp_S(Y)$. The latter follows
directly from (1).

If, in addition, $Y$ is metrizable, then by Corollary
\ref{seccwap}, $(S,Y)$ is RN.
\end{proof}

For every fixed $S$, the class of all RN-approximable compact
$S$-flows is closed under subdirect products. Therefore, using a
well-known method (see for example \cite{Ju,vrbook2}) we obtain that
for every compact $S$-flow $X$ there exists a universal RN-approximable
compactification $u_A: X \to X^A$ which is a topological embedding
iff $\Asp_S(X)$ separates points and closed subsets. Indeed, by
Corollary \ref{compasp} and
Proposition \ref{p-factorRN}.1, it is easy to see that $u_A: X \to
X^A$ is a compactification of $X$  associated to the algebra
$\Asp_S(X)$.

\begin{prop} \label{p-semigr}
For every compact RN-approximable $S$-flow $X$,
its Ellis semigroup $E(X)$, as
an $S$-flow, is also RN-approximable.
\end{prop}
\begin{proof} By the definition, $E(X)$ is an $S$-subflow of $X^X$.
 Hence, $E(X)$ is a subdirect product of
RN-approximable $S$-flows (``$X$ many copies'' of the flow $X$).
\end{proof}

\begin{prop} \label{scattered}
Every scattered compact jointly continuous $S$-flow $X$ is RN.
\end{prop}
\begin{proof}
It is well known that $X$ is scattered iff $C(X)$ is Asplund. Hence
the canonical representation $S \to \Theta(V)_s$,
$X \to B(V^*)_{w^*}$ into an Asplund space $V:=C(X)$ is the desired.
\end{proof}

%%%Let $G$ be a group and  $f: G \to \R$ be a minimal function . Then if
%this function is an Asplund map then it necessarily is almost periodic.
%%%This follows from

Let a semitopological group $G$ act joint continuously on compact
$X$. The following scheme gives some intuitive explanation about
the real place of Asplund functions.

%nd
%\begin{diagram}[height=3.2em,width=4.3em,tight]
% \text{functions:} \hskip 0.5cm & WAP_G(X) & \subset
%& Asp_G(X) & \subset & C(X)\\
% \text{compactifications:} \hskip 0.7cm &  X^W  &\leftarrow & X^{A}
%&\leftarrow &  X^R=X\\
%  \text{representations:} \hskip 1cm
%   & REFL & \subset  & ASP  & \subset &  BAN\\
%\end{diagram}

\vskip 0.1cm

Now define Asplund functions on a semitopological group $G$ via the
universal compactification $u_R: G \to G^R$
(we identify $G$ with $u_R(G)$).
A continuous bounded function $f: G \to \R$ is said to be an {\it Asplund
function} (and write: $f \in \Asp(G)$)
if there exists an Asplund function $F: G^R \to \R$ on the $G$-flow $G^R$
such that $f = F \circ u_R$.
It is equivalent to say that the orbit $fG$ is an Asplund set in the
Banach space $\RUC(G)$. In particular, $\Asp(G) \subset \RUC(G)$.

Define by $\lambda: C(G^R) \to \RUC(G)$
the natural isomorphism of
$G$-algebras induced by the compactification
$u_R: G \to G^R$.

\begin{prop} \label{p-asponG}
Let $G$ be a semitopological group.
\ben
\item
$\Asp(G)$ is a
 $G$-invariant Banach subalgebra of $\RUC(G)$ canonically
$G$-isomorphic to $\Asp_G(G^R)$.
More precisely, $\Asp(G)=\lambda (\Asp_G(G^R))$.
\item
$\WAP(G) \subset \Asp(G)$.
\item
 Denote by $u_A: G \to G^A$ the $G$-compactification
induced by the algebra $\Asp(G)$.
%Then the corresponding action $G \times G^R \to G^R$
%is jointly continuous.
%\item
Then $u_A: G \to G^A$ is the universal RN-approximable jointly continuous
$G$-compactification of $G$. More precisely, for every jointly continuous
$G$-compactification $\a: G \to Y$ with RN-approximable
compact $G$-flow Y, there exists a (necessarily unique) $G$-map
$\psi: G^A \to Y$ such that $\psi \circ u_A=\a$.
\item
$G^A$ is a
right topological monoid naturally isomorphic to the Ellis semigroup
$E(G, G^A)$ and
$u_A$ is a right topological semigroup compactification of $G$.
\een
\end{prop}
\begin{proof} (1)
 Follows by the definition of $\Asp(G)$.

(2) If $f \in \WAP(G)$, then Fact \ref{f-functions}(ii) guarantees that
$fG \subset \RUC(G)$. Then $fG$, being a relatively weakly compact in $\RUC(G)$,
is necessarily an Asplund set (see \cite{Fa}).

%\delete{
%Therefore, $\lambda(WAP_G(G^R)) \subset \lambda(Asp_G(G^R))$. Since
%$\lambda(Asp_G(G^R))=Asp(G)$, it suffices to show that $WAP(G) \subset
%\lambda(WAP_G(G^R))$. By Fact \ref{f-functions}, $WAP(G) \subset RUC(G)$.
%Therefore for every $f \in WAP(G)$ there exists $F \in C(G^R)$ such that
%$\lambda(F)=f$. The orbit $fG$ is relatively weakly compact in
%$RUC(G) \subset C(G)$.
%Since $\lambda: C(G^R) \to RUC(G)$ is a $G$-isomorphism
%of Banach spaces we obtain that $FG$ is also relatively weakly compact in
%$C(G^R)$. Therefore, $F \in WAP_G(G^R)$.
%}

(3) For joint continuity of the action of $G$ on $G^A$, recall that
$\Asp(G)$ is a $G$-invariant subalgebra of $\RUC(G)$. Universality
follows from the fact that
$G^A$ canonically can be identified with $(G^R)^A$ defined for the
jointly continuous compact $G$-flow $G^R$.

(4) Let $i: G \to E(G^A)$ be the natural homomorphism of $G$ into
the Ellis semigroup of the $G$-flow $G^A$. Consider the orbit map
 $\gamma: E(G^A) \to G^A, \gamma(p)=p(u_A(e))$.
Clearly, $\gamma(i(g))=u_A(g)$ for every $g \in G$.
Therefore $\gamma$ is a
morphism between two compactifications $i:G \to E(G^A)$ and $u_A: G
\to G^A$. It suffices to show that $\gamma$
is an isomorphism of these transitive $G$-flows.
By \cite[D.2]{vrbook2} we need
the existence of a morphism of compactifications in the reverse
direction. We can use Proposition \ref{p-semigr} which states that
 $E(G^A)$ is RN-approximable.
By the universality property of $u_A$, there exists
a continuous $G$-map $\nu: G^A \to E(G^A)$ such that
$\nu \circ u_A=i$.
Hence, $\nu$ is the desired morphism between the
compactifications.
\end{proof}

The $G$-algebra $\Asp(G)$ is {\it m-admissible} in the sense of
\cite{BJM} as it follows by Proposition \ref{p-asponG} and
\cite[Theorem 3.1.7]{BJM}.

\begin{thm} \label{matrixasp}
Let $G$ be a semitopological group. The following are equivalent:
\bit

\item [(i)]
$f \in \Asp(G)$.

%\delete{
%(ii) $fG$ is an Asplund set in $RUC(G)$.
%
%(iii) There exist: a $G$-compactification
%$\a: G \to Y$ with a jointly continuous RN-approximable $G$-flow $Y$ and a
%continuous function $F: Y \to \R$ such that $f(g)=F(\a(g))$.
%}

\item [(ii)]
There exist: an Asplund space $V$, a strongly
continuous antihomomorphism $h: G \to \ \Is(V)$, vectors $v \in V$,
and $\psi \in V^*$ such that $f(g)=\langle v, g \psi\rangle $ (that is, $f=m_{v,
\psi}$).

\item [(iii)]
There exist a $G$-compactification $\a: G \to Y$
with a jointly continuous RN $G$-flow $Y$ and a function $F \in
C(Y)=\Asp_G(Y)$ such that $f=F \circ \a$. \item [(iv)] There exist
a $G$-compactification $\a: G \to Y$ with a jointly continuous
$G$-flow $Y$ and a function $F \in \Asp_G(Y)$ such that $f=F \circ
\a$. \eit
\end{thm}
\begin{proof} (i) $\Longrightarrow$ (ii)
Since the $G$-action on a compact space $X=G^R$ is jointly
continuous,  we can use ``strongly continuous version'' of Theorem
\ref{RNrepres} obtaining strongly continuous antihomomorphism $h:
G \to \ \Is(V)$ and a weak$^*$ continuous $\alpha_0 : G^R \to
B(V^*)$ such that $f(g)=\langle v, \alpha_0(g)\rangle $. Then $f=m_{v,\psi}$
where $\psi = \alpha_0(e)$ and $e$ is the identity of $G$.

(ii) $\Longrightarrow$ (iii) Define $Y$ as the weak$^*$ closure of
$cl_{w^*} (G \psi)$ of the orbit of $G \psi$. By Fact
\ref{Teleman}, $Y$ is a jointly continuous $G$-flow. Moreover, $Y$
is an RN $G$-flow (Definition \ref{def-repr}). Then, $F: Y \to
\R$, $F(y)=\langle v,y\rangle $ is the desired function by Corollary
\ref{compasp}.

(iii) $\Longrightarrow$ (iv) is trivial.

(iv) $\Longrightarrow$ (i) By our assumption $FG$ is an Asplund
set in $C(Y)$. The natural $G$-embedding (induced by $\a$) of
Banach algebras $C(Y) \to \RUC(G)$ maps $FG$ onto $fG$. Therefore,
$fG$ is an Asplund set in $\RUC(G)$.
\end{proof}

\begin{prop} \label{p:ap=asp}
Let $G$ be a semitopological group and $X$ be a compact minimal
$G$-flow. Then $AP_G(X)=\WAP_G(X)=\Asp_G(X)$.
\end{prop}
\begin{proof} In general, $AP_G(X) \subset \WAP_G(X) \subset
\Asp_G(X)$ by Lemma \ref{lemma-functions}. So we have only to show
that in our situation $AP_G(X) \supset \Asp_G(X)$ holds. Let $F \in
\Asp_G(X)$. Then by the universality of the canonical $S$-quotient
$u_A: X \to X^A$, we have $F=f \circ u_A$ for some $f \in
\Asp_G(X^A)$. Since $(G,X^A)$ is also minimal, $(G,X^A)$ is
equicontinuous by Theorem \ref{asplundr}. Hence, $f \in
AP_G(X^A)$. Then, clearly $F \in AP_G(X)$. 
\end{proof}

\begin{examples} \label{skew} \ 
\ben
\item
As it was mentioned above, the two-point compactification
%$X=\Z^{+ \infty}_{-\infty}$
of $\Z$, as a cascade, is RN but not wap. Similarly
the two-point compactification
%$\R^{+ \infty}_{-\infty}$
of $\R$, as an $\R$-flow, is RN but not wap.
\item
Define $f: \Z \to \R$ by $f(z)=1$ iff $z$ is a positive integer
and $f(z)=0$ otherwise. Then $f \in \Asp(\Z) \setminus \WAP(\Z)$.

Indeed, $f \in \Asp(\Z)$ by Theorem \ref{matrixasp}
because $f$ comes from
the two-point compactification $Y$ of $\Z$ which is RN $\Z$-flow by (1).
On the other hand, $f \notin \WAP(\Z)$ because $f$ does not satisfy
Grothendieck's DLP (see Fact \ref{f-DLP}). Choose $s_n=n$ and $x_m=-m$.
Then $\lim_m \lim_n f(n-m) = 1 \neq 0 = \lim_n \lim_m f(n-m)$.

\item As in (2), it is easy to show that
$f\in \Asp(\R) \setminus \WAP(\R)$ for the functions $f(x) = \frac{x}{1+|x|}$
%zaza270413 replacing arcsin
and $f(x)=arctg x$.

\item The cascade $(\Z, [0,1])$ generated by the $f(x)=x^2$ map is
not wap. Indeed, it contains, as a subflow, the two-point
compactification of $Z$. Take, for example, the $\Z$-orbit of the
point $x=\frac{1}{2}$. Together with the endpoints $\{0\}$ and
$\{1\}$, we get the closure of this orbit.

%\delete{
%\item
%Let $\R_+=(0,\infty)$ be the multiplicative group of all positive reals.
% Consider the natural action $\R_+ \times \R \to \R, (t,r) \mapsto tr$.
%Define $f: \R \to \R$ by $f(x)=x$ for all $|x| < 1$ and $f(x)=1$ otherwise.
%Then $f \in Asp_{\R_+}(\R) \setminus WAP_{\R_+}(\R)$. Moreover,
%it follows directly by the definition of Asplund functions that
%every function $f \in C_0(\R)$ vanishing at infinity
%belongs to $Asp_{\R_+}(\R)$. Indeed, one can check that
%the set $\Q$ of all rationals is dense in
%the semimetric space $(\R, \rho_f)$. At the same time there is no
%function from $WAP_{\R_+}(\R)$ which separates the point $0$ and the subset
%$\R \setminus (-\e,\e)$ for a fixed $\e >0$.
%Indeed, assuming the contrary, let $f(0)=0$ and
%$f(x)=1$ for every $|x| \geq \e$. Consider the sequences $g_m=m \in \R_+$ and
%$x_n=\frac{1}{n} \in \R$. Then $\lim_m \lim_n f(g_m x_n)=0$ and
%$\lim_n \lim_m f(g_m x_n)=1$. Therefore $f \notin WAP_{\R_+}$ because
% DLP is not satisfied (see Fact \ref{f-DLP}).
%It follows that $(\R_+, \R)$
%(and even any ``multiply by scalars flow''
%$(\R_+, V)$ for every normed space $V$) is RN but not wap.
%}

\item

Let $X$ be a minimal compact jointly continuous $G$-flow which is
not equicontinuous. Then $C(X) \neq \Asp_G(X)$ ($X$ is not
RN-approximable) and $\RUC(G) \neq \Asp(G)$.

%30/10/02
Indeed, by Theorem \ref{asplundr}, $X$ is not RN-approximable.
Theorem \ref{compasp} guarantees that $C(X) \neq \Asp_G(X)$.
Another proof of the same fact follows from the equality
$AP_G(X)=\Asp_G(X)$ (Proposition \ref{p:ap=asp}).

Now we check that $\RUC(G) \neq \Asp(G)$. Fix $f \in C(X) \setminus
\Asp_G(X)$ and a point $z \in X$. Since $z$ is a point of
transitivity of $X$, there exists a continuous onto $G$-map $q:
G^R \to X$ such that $q(u_R(g))=gz$ for every $g \in G$. Define
$F: G^R \to \R$ as the composition $f \circ q$. Then $F \notin
\Asp_G(G^R)$ by Proposition \ref{p-factorRN}. Thus the restriction
$F|_G(g)=f(gz)$ of $F$ on $G$ satisfies $F|_G \in \RUC(G) \setminus
\Asp(G)$.

As a concrete example consider the cascade on the
two-dimensional torus $\T^2=(\R/\Z)^2$
generated by the selfhomeomorphism (see \cite[Example 5.1.7]{BJMo})
$$\sigma_{\theta}: \T^2 \to \T^2, \hskip 0.2cm
\sigma_{\theta} ([a],[b]) = ([a+b], [a+\theta]) $$ where $\theta$
is a given irrational number. Then the corresponding flow $(\Z,
\T^2, \pi_{\theta})$ is minimal but not equicontinuous. The
minimality one can check by results of Furstenberg \cite{Fu}. In
particular, the cascade $(\Z, \T^2, \pi_{\theta})$ is not RN. Another corollary: $\Asp(\Z) \neq \RUC(\Z)=C(\Z)$ and $\Z^R$ is
not RN-approximable. \een
\end{examples}

\begin{remarks} \label{r-Ryll} \ 
\bit
\item [(i)]
A result from \cite{MPU} states that a topological group $G$ is precompact
iff $\WAP(G)=\RUC(G)$,
previously  obtained in \cite{AG} for monothetic groups.
Is it true the same assuming $\Asp(G)=\RUC(G)$ ?
\item [(ii)]
Namioka and Phelps \cite{NP} proved a generalized Ryll-Nardzewski
fixed-point theorem for $S$-flows which are weak star compact
convex subsets in the dual of an Asplund space. Hence, this
situation is a particular case of RN flows. It is interesting to
analyze possible applications for amenability context as well as
for decomposition theorems.
\eit
\end{remarks}

\section{Kadec property: ``when does weak imply strong ?''}

\vskip 0.1cm

 Recall that a Banach space $V$ has the {\it Kadec property} if
the weak and norm topologies coincide on the unit (or some other)
sphere of $V$. Let us say that a subset $X$ of a locally convex
space (l.c.s.) $(V,\tau)$ is a {\it Kadec subset} ({\it light
subset} in \cite{Meop}) if the weak topology coincides with the
strong topology. {\it Light} linear subgroups $G \leq Aut(V)$
(with respect to the weak and strong {\it operator topologies})
can be defined Analogously. Clearly, if $G$ is {\it orbitwise
Kadec} on $V$ that is,
%$\tau_w|_X=\tau|_X$
all orbits $Gv$ are light in $V$, then $G$ is necessarily light.
The simplest examples are the spheres
(orbits of the unitary group $\Is(H)$) in Hilbert spaces $H$.

The following results show that linear actions frequently are
``orbitwise Kadec''.

\begin{thm} \label{kadec}
Let a subgroup $G \leq Aut(V)$ be
equicontinuous, $X$ be a bounded,
$(weak, \mu)$-fragmented
$G$-invariant subset of an l.c.s. $V$ with the natural uniformity
$\mu$. Then every, not necessarily closed, quasiminimal
$G$-subspace (e.g., the orbits) $Y$ of $X$ is a Kadec subset.
\end{thm}
\begin{proof} The equicontinuity of the subgroup $G\leq Aut(V) $
implies that the action of $G$ on the bounded subspace $X \subset V$ is
uniformly $\mu$-equicontinuous
with respect to the natural uniformity $\mu$ on $V$. Since $X$ is
$(weak, \mu)$-fragmented we get that in fact
the $G$-flow $X$ is (weak, $\mu|_X$)-equifragmented.
Therefore we can apply Theorem \ref{genquasi}.
\end{proof}

%\begin{defin}
We say that an l.c.s. $V$ is {\it boundedly fragmented} (write: $V
\in BF$) if every bounded subset $X \subset V$ is
(weak, $\mu$)-fragmented, where $\mu$, as above, is the natural
uniformity of $V$.
%\end{defin}

\begin{corol} \label{PCP}
Let $V \in BF$. Then every equicontinuous $G \leq Aut(V)$ is a
light subgroup and every orbit $Gv$ is a light subset in $V$.
\end{corol}

The class BF is large (see the relevant references in
\cite{Meop}) and includes among others:
%\nt {\bf Examples of (BF)}:
Banach spaces with PCP (point of continuity property),
%(by [Jayne-Rogers 1985]).
semireflexive l.c.s., Frechet spaces with the Radon-Nikod\'ym
Property.
%(by [Egghe 1977]).

%The following result has been obtained also in \cite{}.

\begin{corol} \label{old}
\cite{Meop} Let $V$ be a Banach space with PCP. Then any bounded
subgroup $G$ of $Aut(V)$ (e.g., $\Is(V)$) is light.
\end{corol}

Now, combining Corollary \ref{old} and Theorem \ref{wapeb},
one can obtain a soft geometric proof of Fact \ref{lawson}
using representations of dynamical systems on reflexive spaces. 

\vskip 0.1cm
 \nt {\it Ellis-Lawson's Joint Continuity
Theorem: Let $G$ be a subgroup of a  compact
semitopological monoid $S$. Suppose that $(S,X)$ is a
semitopological flow with compact $X$. Then the action
$G\times X \to X$ is jointly continuous and $G$ is a topological
group.}
\begin{proof} $(S,X)$ is wap by Fact \ref{EN}.
Therefore by Theorem \ref{wapeb}
there exists an approximating family $(h_i, \alpha_i)$ of
reflexive $V_i$-representations. It suffices to prove the theorem
for the canonical wap $\Theta(V)^{opp}$-flow
$B^*$. Let $G$ be a
subgroup of $\Theta(V)^{opp}$. Then by Corollary \ref{old}  the
strong operator topology on $G$ coincides with the weak topology.
In particular, $G$ is a topological group. Moreover, by Fact
\ref{Teleman} the action of $G$ on $(B^*,w^*)$ is jointly
continuous.
\end{proof}

%\nt {\bf  Ellis theorem}. Every locally compact semitopological
%group $G$ is a topological group.

%\delete{
%\nt Proof: locally compact semitopological group $G$ is wap
%(because the Alexandrov compactification $G^*$ is semitopological).
%Then by by the representation theorem of Shtern
%$G^*$ can be embedded into $\Theta(V)_w$ for some reflexive $V$.
%Hence $G$ is embedded into $Is(V)_w$.
%But by Corollary \ref{old} $Is(V)_w=Is(V)_s$. Hence, $G$ is a topological
%group.
%}
%

%(II) Lawson j. c. thm (a weaker version)
%$h: G \to Aut(V), \quad {\tilde v}: S \to V$

Recall that $X$ is a  {\it Namioka space} if for every compact
space $Y$ and a separately continuous map $\gamma: Y \times X \to
\R$ there exists a dense subset $P \subset X$ such that $\gamma$
is jointly continuous at every $(y,p) \in Y \times P$.
A topological space is said to be
\v{Chech}-{\it complete} if it can be
represented as a $G_{\delta}$-subset of a compact space.
Every \v{Chech}-complete  (e.g., locally compact or Polish) space is a
Namioka space.

\begin{prop} \label{chech}
Let $G$ be a semitopological group, $X$ be a semitopological $G$-flow
and $f \in C(X)$. Suppose that $(cl_w(fG), weak)$ be
a {\it Namioka space}. Then $fG$ is light in $C(X)$.
\end{prop}
\begin{proof} As in the proof of Lemma \ref{namiokamy},
it is easy to show that if a bounded subset of a Banach space is a
Namioka space under the weak topology then it is (weak-norm)-fragmented.
Therefore we can complete the proof by Theorem \ref{kadec}.
\end{proof}

\begin{thm} \label{c-waplight}
Let $G$ be a semitopological group. Then for every semitopological
$G$-flow $X$ and every $f \in \WAP(X)$ the pointwise and norm
topologies coincide on the orbit $fG$. In particular, $\WAP_G(X)
\subset \RUC_G(X)$.
\end{thm}
\begin{proof}
Proposition \ref{chech} guarantees that $fG$ is (weak, norm)-Kadec. On
the other hand, the weak and pointwise topologies coincide on the
weak compact set $cl_w(fG)$.
\end{proof}

%This implies in particular an easy proof of the  following well known fact
%$WAP(G) \subset RUC(G)$.

\vskip 0.1cm  Now we turn to the weak$^*$ version of the lightness
concept. Let $(V,\tau)$ be an l.c.s. with its strong dual $(V^*,
\tau^*)$. Denote by $\mu^*$ the corresponding uniformity on $V^*$.
Let's say that a subset $A$ of $V^*$ is {\it weak$^*$ light} if
weak$^*$ and strong topologies coincide on $A$. If $G$ is a
subgroup of $Aut(V^*)$, then the weak$^*$ (resp., strong$^*$)
topology on $G$ is the weakest topology which makes all orbit maps
$\{{\tilde \psi}: G \to V^*: \hskip 0.1cm \psi \in V^* \} $
weak$^*$ (resp., strong) continuous.

Following \cite{Mefr} we say that an l.c.s. $V$ is a {\it
Namioka-Phelps space} ($V \in NP$) if every equicontinuous subset
$X \subset V^*$ is $(w^*, \mu^*)$-fragmented.
%\end{defin}
The class NP is closed under subspaces, products and l.c. sums and includes:
Asplund Banach spaces, semireflexive l.c.s. and Nuclear l.c.s.

\begin{thm} \label{NP}
Suppose that  $V$ is an NP space, $G\leq Aut(V)$ is an
equicontinuous subgroup, and  $X \subset V^*$ is an equicontinuous
$G$-invariant subset.
\bit
\item [(i)] If $(X, w^*)$ is a quasiminimal (e.g.,
1-orbit) G-subset, then $X$ is weak$^*$ light.
%(that is, $w^* |_X=\tau^* |_X$).
\item
[(ii)] The weak$^*$ and strong$^*$ operator topologies coincide on
$G$. \eit
\end{thm}
\begin{proof} (i) The strong topology on the dual space $(V^*, \mu^*)$
is the topology of bounded convergence.
Since $G$ is an equicontinuous subgroup of
$Aut(V)$, it is easy to show that the dual action of $G$ on $V^*$ is also
equicontinuous. On the other hand,
$X \subset V^*$ is $(w^*, \mu^*)$-fragmented as it follows
by the definition of NP spaces. We obtain in fact that
$G$-flow $X$ is (weak$^*$, $\mu^*|_X$)-equifragmented.
Now use once again Theorem \ref{genquasi}.

(ii) Directly follows from (i) because every $G$-orbit in $V^*$
is an equicontinuous subset.
\end{proof}

The last result is useful in the context of continuity of {\it dual actions}
(for more information see \cite{Mefr} and the references there).
More precisely,
let $V$ be an l.c.s. and $h: G \to Aut(V)$ be a homomorphism such that
$h(G)$ is an equicontinuous subgroup of $Aut(V)$ and the action
$G \times V \to V$
is jointly continuous. Then we can ask: is the dual action
$$\pi^*: G \times (V^*, \mu^*) \to (V^*, \mu^*), \hskip 0.1cm (gf)(v)=f(g^{-1}v)$$
also jointly continuous ?

Since $h(G)$ is equicontinuous, clearly $(G, V^*)$
is equicontinuous with respect to the dual action $\pi^*$. Therefore
it is equivalent to ask if the orbit maps
${\tilde f}: G \to (V^*, \mu^*)$ are continuous for all $f \in V^*$.
Since ${\tilde f}: G \to (V^*, w^*)$ is continuous,
it suffices to show that the orbits $Gf$ are (weak$^*$, strong)-Kadec subsets
of $V^*$.
This fact follows directly from Theorem \ref{NP} provided that $V \in NP$.
Hence we obtain the following result.

\begin{corol} \label{c-dual}
Let $V \in NP$ (e.g., Asplund Banach space)
and $\pi: G \times V \to V$ be a linear
jointly continuous equicontinuous action. Then the dual action
$\pi^*: G \times V^* \to V^*$ is also jointly continuous.
\end{corol}

\begin{remark} \label{r-dual}
Corollary \ref{c-dual} can be derived also from
\cite[section 6]{Mefr}.
If $V$ is an Asplund Banach space then we can drop
the condition about equicontinuity
(in fact, boundedness) as it follows by
\cite[Corollary 6.9]{Mefr}.
\end{remark}

\begin{prop} \label{right}
Let $V$ be an Asplund Banach space, $G$ be a semitopological group,
and $h: G \to Aut(V)$ be a
bounded weakly continuous antihomomorphism.
Assume that $v \in V$ and $\psi \in V^*$. 
Then the corresponding matrix coefficient
$m_{v,\psi}: G \to \R$ is left uniformly continuous. Moreover, if the vector
$v$ is norm-continuous, then $m_{v,\psi}$, in addition,
is right uniformly continuous.
\end{prop}
\begin{proof}
The antihomomorphism $h$ sends $G$ into a norm bounded subgroup of
$Aut(V)$. Therefore by Fact \ref{introversion} it suffices to show
that $\psi$ is a norm $G$-continuous vector. Since $h: G \to
Aut(V)$ is weak continuous, the orbit map $\tilde{\psi}: G \to
V^*$ is weak star continuous. Since $V$ is Asplund (and hence NP),
Theorem \ref{NP} implies that the weak star and norm topologies
coincide on the orbit $G \psi$. Then $\tilde{\psi}: G \to V^*$ is
even norm continuous.
\end{proof}

%\delete{
%For every $\varepsilon >0$ there exists a nbd $O(e)$ such that
%$t\psi - \psi$ is $\varepsilon$-small. Then
%$gt\psi - g\psi$ is $\|G\|\cdot \varepsilon$-small for every $g \in G$. Then
%$<v,gt\psi>-<v,g\psi>= m(gt)-m(g)$
%is $\|v\| \cdot \|G\|\cdot \varepsilon$-small. Therefore, $m$ is
%right uniformly continuous.
%In the case that the orbit map $\tilde{v}: G \to V$ is norm continuous ,
%we can choose $O(e)$ such that $<vt-v, g\psi>=m(tg)-m(g)$
%is sufficiently small for every $g \in O(e)$.
%}

%Combining Theorems \ref{right} and \ref{matrixrn} we obtain

\begin{corol} \label{forasplf}
For every semitopological group $G$  $$\WAP(G) \subset \Asp(G)
\subset \LUC(G) \cap \RUC(G)$$
 holds.
\end{corol}
\begin{proof}
The inclusion $\WAP(G) \subset \Asp(G)$
is a part of Proposition \ref{p-asponG}.
Let $f \in \Asp(G)$. By Theorem \ref{matrixasp} the function $f$ coincides
with a matrix coefficient $m_{v, \psi}$ for a suitable strongly continuous
antihomomorphism $h: G \to \Is(V)_s$. Now we can apply
Proposition \ref{right} to $f=m_{v, \psi}$. 
\end{proof}

%\section{Other applications}

%%%%%%%%%%%%%%%%%%     REFERENCES

\bibliographystyle{amsplain}

\begin{thebibliography}{10}


\bibitem{AAB}  E. Akin, J. Auslander and K. Berg,
\textit{Almost equicontinuity and the enveloping semigroup},
Contemp. Math., \textbf{215} (1998), 75-81.

\bibitem{AG}  E. Akin and E. Glasner,
{\it Residual properties and almost equicontinuity},
J. Anal. Math., {\bf 84} (2001), 243-286.

\bibitem{AL} D. Amir and J. Lindenstrauss,
{\it The structure of weakly compact subsets in Banach spaces},
Ann. of Math., {\bf 88} (1968), 35-46.

%\bibitem{AMM} I. Aharoni, B. Maurey and B.S. Mityagin,
%{\it uniform embedding of metric spaces and of Banach
%spaces into Hilbert spaces},
%Israel J. Math., {\bf 52} (1985), 251-265.

\bibitem{Au}  J. Auslander,
{\it Minimal flows and their extensions},
{\bf 153}, North-Holland Math. St., 1988.

\bibitem{Arh} A.V. Arhangel'ski\u\i,
{\it On some topological spaces that occur in Functional Analysis,}
Russ. Math. Surveys,
{\bf 31:5} (1976), 17-32.

\bibitem{BL}
Y. Benyamini and J. Lindenstrauss,
{Geometric Nonlinear Functional Analysis},
vol. 1, Amer. Math. Soc., Colloquium Publications {\bf
48}, Providence, Rhode Island, 2000.

%\delete{
%\bibitem{Ba} W. Banasczyk,
%{\it Additive subgroups of topological vector spaces},
%Lecture Notes in Math., {\bf 1466} (1991), Springer-Verlag.
%}


\bibitem{BJMo}  J.F. Berglund, H.D. Junghenn and P. Milnes,
{\it Compact right topological semigroups and
generalizations of almost periodicity}, Lecture Notes in Math.,
{\bf 663} (1978), Springer-Verlag.

\bibitem{BJM}  J.F. Berglund, H.D. Junghenn and P. Milnes,
{\it Analysis on Semigroups}, Wiley, New York, 1989.

\bibitem{Bour} R.D. Bourgin,
{\it Geometric aspects of of convex sets with the Radon-Nikod\'ym
property}, Lecture Notes in Math. {\bf 993}, Springer-Verlag,
1983.

\bibitem{Bo} A. Bouziad,
{\it Continuity of separately continuous group actions in $p$-spaces},
 Topology Appl. {\bf 71} (1996), no. 2, 119--124.

\bibitem{Bu} R.B. Burckel,
{\it Weakly almost periodic functions on semigroups},
Gordon and Breach Science Publishers, New York-London-Paris, 1970.

\bibitem{DFJP}
W.J. Davis, T. Figiel, W.B. Johnson and A. Pelczynski, {\em
Factoring weakly compact operators\/}, {\it J. of Funct. Anal.},
{\bf 17} (1974), 311-327.

\bibitem{DGZ} R. Deville, G. Godefroy and V. Zizler,
{\it Smoothness and renormings in Banach spaces},
Pitman Monographs and Surveys in Pure and Applied Mathematics, {\bf 64},
Longman Scientific Technical, 1993.

\bibitem{Do} T. Downarowicz,
{\it Weakly almost periodic flows and hidden eigenvalues},
Contemp. Math., {\bf 215} (1998), 101-120.

\bibitem{Eb} W.F. Eberlein,
{\it Abstract ergodic theorems and weak almost periodic functions},
Trans. AMS, {\bf 67} (1949), 217-240.

%\delete{
%\bibitem{Eg} L. Egghe, {\it On the Radon-Nikod\'ym property, and related
%topics in locally convex spaces}, Lecture Notes in Math., {\bf 645}
%(1977), 77-90.
%}

\bibitem{EN} R. Ellis and M. Nerurkar,
{\it Weakly almost periodic flows},
Trans. AMS, {\bf 313} (1989), 103-119.

%\delete{
%\bibitem{En} P. Enflo,  {\it On a problem of Smirnov}, Ark. Math.,
%{\bf 8} (1969), 107-109.
%}
%

\bibitem{Fa}  M. Fabian,
{\it Gateaux differentiability of
convex functions and topology.
Weak Asplund spaces}, Canadian Mathematical Society
Series of Monographs and Advanced Texts,
Wiley, New York, 1997.

%
%\delete{
%\bibitem{Fu}  H. Furstenberg,
%{\it The structure of distal flows},
%Amer. J. Math., {\bf 85} (1963), 477-515.
%}
%
%\delete{
%\bibitem{FG} H. Furstenberg, E. Glasner,
%{\it Robert Ellis and the algebra of dynamical systems},
%Contemp. Math., {\bf 215} (1998), 197-203.
%}

\bibitem{Fu}  H. Furstenberg,
{\it Strict ergodicity and transformation of the torus},
Amer. J. Math., {\bf 83} (1961), 573-601.

\bibitem{Gl} E. Glasner,
{\it $M$-Dynamical systems}, preprint, 1998.

\bibitem{GW1} E. Glasner and B. Weiss,
{\it Sensitive dependence on initial conditions},
Nonlinearity {\bf 6} (1993), 1067-1075.

\bibitem{GW2} E. Glasner and B. Weiss,
{\it Locally equicontinuous dynamical systems},
Colloquium Mathematicum, part 2, {\bf 84/85} (2000), 345-361.

\bibitem{Gr} A. Grothendieck,
{\it Crit\'eres de compacit\'e dans les espaces
functionelles g\'en\'eraux},
Amer. J. Math., {\bf 74} (1952), 168-186.

%\bibitem{Gu} I. Guran,
%{\it On topological groups close to being Lindelof},
%Soviet Math. Dokl., {\bf 23} (1981), 173-175.

\bibitem{HT}
G. Hansel and J.P. Troallic,
{\it Extension properties
of WS-groups}, Semigroup Forum, {\bf 45:1} (1992), 63--70.

\bibitem{He1} D. Helmer,
\emph{Joint continuity of affine semigroup actions}, Semigroup
Forum, {\bf 21} (1980), 153-165.

\bibitem{He2} D. Helmer,
{\it Continuity of semigroup actions},
Semigroup Forum, {\bf 23} (1981), 153-188.

%\delete{
%\bibitem{HC} W. Herer, J.P.R. Christensen,
%{\it On the existence of pathological submeasures and the construction
%of exotic topological groups}, Math. Ann., {\bf 213} (1975), 203-210.
%
%%\bibitem{Ka} M. Katetov, {\it On
%%real-valued functions in topological spaces},
%%Fund. Math. {\bf 38}, 85-91, and correction {\bf 40} (1953), 203-205.
%}

\bibitem{Ke} P. Kenderov,
{\it Dense strong continuity of pointwise continuous mappings},
Pacific. J. Math., {\bf 89:1} (1980), 111-130.

\bibitem{KM} J.L. Krivine and B. Maurey,
{\it Espaces de Banach stables},
Israel J. Math. {\bf 39:4} (1981), 273--295.

\bibitem{JP} A. Jalilian and  M.A.  Pourabdollah,
{\it Transformation semigroup compactifications and norm
continuity of weakly almost periodic functions},
Proc. Indian Acad. Sci. Math. Sci. {\bf 110:1} (2000), 55--59.
%\delete{
%\bibitem{JNR} J.E. Jayne, I. Namioka and C.A. Rogers,
%{\it Topological properties of Banach spaces}, Proc.
%London Math. Soc. (3), {\bf 66}, 1993, 651-672.
%}

\bibitem {JOPV}
J.E. Jayne, J. Orihuela, A.J. Pallares and G. Vera,
{\it $\sigma$-fragmentability of multivalued maps and
selection theorems},
J. Funct. Anal. {\bf 117} (1993), no. 2, 243--273.

\bibitem{JR} J.E. Jayne and C.A. Rogers,
{\it Borel selectors for upper semicontinuous set-valued maps},
Acta Math., {\bf 155}, (1985), 41-79.

\bibitem{Ju}  H.D. Junghenn,
{\it Almost periodic compactifications of product flows},
Semigroup Forum, {\bf 58} (1999), 296-312.

\bibitem{La1}  J.D. Lawson,
{\it Joint continuity in semitopological semigroups},
Illinois J. Math., {\bf 18} (1974), 275-~285.

\bibitem{La2}  J.D. Lawson,
{\it Points of continuity for semigroup actions},
Trans. AMS {\bf 284} (1984), 183-202.

\bibitem{LG1} K. de Leeuw and I. Glicksberg,
{\it Applications of almost periodic compactifications},
Acta Math., {\bf 105} (1961), 63-97.

%\delete{
%\bibitem{LG2} K. de Leeuw and I. Glicksberg,
%{\it The decomposition of certain group
%representations}, Journal D'Analyse, {\bf 15} (1965), 135-192.
%}
%%\bibitem{Li} J. Lindenstrauss,
%{\it Uniform embeddings, homeomorphisms and
%quotient maps between Banach spaces (a short survey)}
%Preprint, 1997.

%\bibitem{Li} J. Lindenstrauss,
%{\it Weakly compact sets : their topological properties and
%the Banach space they generate},
%Annals of Math. Studies {\bf 69}, Princeton Univ. Press, 1972,
%235-273.

%delete{
%\bibitem{Ma} P. Mankiewicz, {\it On differentiability
%of Lipschitz mappings in Fr\'echet spaces},
%Studia Math.,  {\bf 45} (1973), 15-29.
%
%\bibitem{May} H. Maynard, {\it A geometric characterization of
%Banach spaces possessing the Radon-Nikod\'ym-Property},
%Trans. AMS, {\bf 185} (1973), 493-500.
%}


\bibitem{LNO}
A. Lima, O. Nygaard, E. Oja,
\emph{Isometric factorization of weakly compact operators and the approximation property,}
Israel J. Math., \textbf{119} (2000), 325-348.

\bibitem{Mesm} M. Megrelishvili,
{\it A Tychonoff $G$-space not admitting a compact $G$-extension or a
$G$-linearization}, Russian Math. Surv., {\bf 43:2} (1988), 177-178.

%\delete{
%\bibitem{Memin}  M.G. Megrelishvili, {\it Group
%representations and construction of minimal topological groups},
%Topology Appl., {\bf 62}, 1995, 1-19.
%}

\bibitem{Mefr} M. Megrelishvili,
{\it Fragmentability and continuity of semigroup actions},
Semigroup Forum, {\bf 57} (1998), 101-126.

%\delete{
%\bibitem{MeF}  M.G. Megrelishvili, {\it Free Topological G-Groups},
%New Zealand Journal of Mathematics, vol. 25, no. 1, 1996, 59-72.
%
%\bibitem{Meeb}  M.G. Megrelishvili,
%{\it Eberlein groups and compact semitopological semigroups},
%Bar-Ilan University preprint, 1998.
%}

%\bibitem{Meopold}  M.G. Megrelishvili,
%{\it Operator topologies and reflexive
%representability of groups},
%Bar-Ilan University preprint, 1998.

\bibitem{Merup} M. Megrelishvili,
{\it Every semitopological semigroup compactification of the group
$H_+[0,1]$ is trivial},
Semigroup Forum, {\bf 63:3} (2001), 357-370.

\bibitem{Meop}  M. Megrelishvili,
{\it Operator topologies and reflexive
representability},
In: ``Nuclear Groups and Lie Groups'',
Research and Exposition in Math. Series,
{\bf 24} (2001), Heldermann-Verlag, 197-208.

\bibitem{Meist}  M. Megrelishvili,
{\it Reflexively but not unitarily representable topological groups},
Topology Proceedings, {\bf 25} (2002), 615-625.

\bibitem{MeSc} M. Megrelishvili and T. Scarr,
{\it Constructing Tychonoff G-Spaces Which Are Not $G$-Tychonoff},
Topology and its Applications, {\bf 86:1} (1998), 69-81.

\bibitem{MPU}  M. Megrelishvili, V. Pestov and V. Uspenskij,
{\it A note on the precompactness of weakly almost periodic groups}
In: ``Nuclear Groups and Lie Groups'',
Research and Exposition in Math. Series,
{\bf 24} (2001), Heldermann-Verlag, 209-216.

\bibitem{Mi} P. Milnes,
{Compactifications of semitopological semigroups},
J. Austral. Math. Soc. {\bf 15} (1973), 488--503.

%\delete{
%\bibitem{Mo}  R.T. Moore, {\it Measurable, continuous and
%smooth vectors for semigroup and group representations},
%Mem. Amer. Math. Soc., {\bf 78} (1968).
%}

\bibitem{Na1} I. Namioka,
{\it Separate continuity and joint
continuity}, Pacific. J. Math., {\bf 51} (1974), 515-531.

\bibitem{Na2} I. Namioka,
{\it Radon-Nikod\'ym compact spaces and fragmentability},
Mathematika, {\bf 34} (1987), 258-281.

\bibitem{NP}  I. Namioka and R.R. Phelps,
{\it Banach spaces which are Asplund spaces},
Duke Math. J., {\bf 42} (1975),  735-750.

\bibitem{NW} I. Namioka and R.F. Wheeler,
{\it Gulko's proof of the Amir-Lindenstrauss theorem},
Contemp. Math., {\bf 52} (1986), 113-120.

\bibitem{Pe1}  V.G. Pestov,
{\it Topological groups: where to from here?}
Topology Proceedings,
{\bf 24} (1999), 421-502. http://arXiv.org/abs/math.GN/9910144.

\bibitem{Re} O.I. Reynov,
{\it On a class of Hausdorff compacts and GSG Banach spaces},
Studia Math., {\bf 71} (1981), 113-126.

\bibitem{Ro1} H.P. Rosenthal,
{\it The heredity problem for weakly compactly generated Banach spaces},
Compositio Math. {\bf 28} (1974), 83--111.

\bibitem{Ru}  W. Ruppert,
{\it Compact semitopological semigroups: An intrinsic theory},
Lecture Notes in Math., {\bf 1079} (1984), Springer-Verlag.

%\bibitem{Ru2} W. Ruppert,  {\it
%On signed $a$-adic expansions and weakly almost
%periodic functions}, Proc. London Math. Soc. (3), {\bf 63} (1991), 620-656.

%\bibitem{Sho}
%I.J. Shoenberg, Metric spaces and positive definite
%functions, Trans. AMS {\bf 44}, 1938, 522-536.

\bibitem{Sh}  A. Shtern,
{\it Compact semitopological semigroups and
reflexive representability of topological groups},
Russian J. of Math. Physics, {\bf 2} (1994), 131-132.

\bibitem{St} Ch. Stegall,
{\it The Radon-Nikod\'ym property and conjugate Banach spaces} II,
Trans. AMS, {\bf 264} (1981), 507-519.

%\bibitem{St} C. Stegall, {\it The duality between Asplund spaces and
%spaces with the Radon-Nikod\'ym property}, Israel J. Math., {\bf 29} (1978),
%408-412.
%\delete{
%\bibitem{Sz} W. Szlenk, {\it The non-existence of a separable reflexive
%Banach space universal for all separable reflexive Banach spaces},
%Studia Math., {\bf 30} (1968), 53-61.
%}

\bibitem{Te} S. Teleman,
{\it Sur la repr\'esentation lin\'eare des groupes topologiques},
Ann. Sci. Ecole Norm. Sup., {\bf 74} (1957), 319-339.

%\bibitem{Tr} J.P. Troallic, {\it Semigroupes semitopologiques et
%presqueperiodicite}, Lecture Notes in Math., {\bf 998} (1981),
%Springer-Verlag, 239-251.

\bibitem{Tr} J.P. Troallic,
{\it Espaces Fonctionnels et theoremes de I.
Namioka}, Bull. Soc. Math. France, {\bf 107} (1979), 127-137.

%\delete{
%\bibitem{Us1}  V.V. Uspenskij,
%{\it The Roelcke compactification of unitary groups},
%in: Abelian groups, module theory, and topology,
%Proceedings in honor of Adalberto Orsatti's 60th birthday
%(D.  Dikranjan, L. Salce, eds.),
%Lecture notes in pure and applied mathematics,
%Marcel Dekker, New York e.a., {\bf 201} (1998), 411-419.
%
%
%%\bibitem{Us2}  V.V. Uspenskij,
%%{\it The Roelcke compactification of groups of homeomorphisms},
%%Topology and its Applications,
%%{\bf 111} (2001), 195-205.
%}

\bibitem{Usold}  V.V. Uspenskij,
{\it Universal topological groups with a countable basis},
Functional Anal. Prilojen. (in Russian), {\bf 20:2} (1986), 86-87.

%\delete{
%
%[Us2] V.V. Uspenskii, On subgroups of minimal topological groups,
%Preprint, 1998, 1-29.
%\item
%
%}

%\bibitem{Usun}  V.V. Uspenskij,
%{\it On universal minimal compact $G$-spaces},
%Proceedings of the 2000 Topology and Dynamics Conference
%(San Antonio, TX), Topology Proceedings, {\bf 25} (2000), 301--308.

\bibitem{Ve}  W. A. Veech,
{\it A fixed point theorem free approach to weak almost
periodicity}, Trans. AMS, {\bf 177} (1973), 353-362.

%\bibitem{vrbook1} J. de Vries, {\it Topological transformation
%groups}, Math. Centre Tract {\bf 65}, Mathematisch Centrum,
%Amsterdam, 1975.

%\bibitem{vrexist}
%J.de Vries, {\it On the existence of $G$-compactifications}, Bull.
%Ac. Polon. Sci. Ser. Math., {\bf 26} (1978), 275-280.

\bibitem{vr-embed} J. de Vries, {\it Equivariant embeddings of
$G$-spaces}, in: J. Novak (ed.), {\it General Topology and its
Relations to Modern Analysis and Algebra IV}, Part B, Prague,
1977, 485-493.

\bibitem{vrbook2} J. de Vries,
{\it Elements of Topological Dynamics},
Kluwer Academic Publishers, Dordrect-Boston-London, 1993.

\end{thebibliography}

\end{document}